\documentclass[10pt]{article}
\usepackage{amsmath,mathrsfs,dsfont,amssymb,epsfig,color,graphicx,amsthm}
\usepackage[latin1]{inputenc}
\usepackage{dsfont}
\newtheorem{theorem}{Theorem}
\newtheorem{proposition}[theorem]{Proposition}
\newtheorem{corollary}[theorem]{Corollary}
\newtheorem{lemma}[theorem]{Lemma}
\newtheorem{definition}{Definition}

\newcommand{\cqd}{\nopagebreak\hfill\fbox{ }}
\let\a=\alpha

\let\L=\Lambda

\newcommand{\R}{\mathbb{R}}

\newcommand{\Z}{\mathbb{Z}}
\newcommand{\N}{\mathbb{N}}

\newcommand{\rar}{\rightarrow}
\newcommand{\LL}{\mathcal{L}}

\newcommand{\TT}{\mathcal{T}_{s,A}}
\newcommand{\LLW}{\mathcal{L}_{A}}
\newcommand{\DLLW}{\mathcal{L}^*_{A}}
\newcommand{\evc}{\psi_A}
\newcommand{\evl}{\lambda_A}
\newcommand{\II}{\int_M}
\newcommand{\BB}{\mathcal{B}}

\begin{document}
\setcounter{section}{-1}

\title{On the general One-Dimensional $XY$ Model: positive and zero temperature, selection and non-selection}
\author{A. T.  Baraviera(*), L. M. Cioletti(**),  A. O. Lopes (*),\\
J. Mohr (*) and R. R. Souza(*)} 
\date{
(*)\,Instituto de Matem\'atica, UFRGS - Porto Alegre, Brasil and (**) Departamento de Matematica, UNB, Brasilia, Brasil.}

\date{\today}
\maketitle

\begin{abstract}

We consider $(M,d)$ a connected and compact manifold and
we denote by $\mathcal{B}_i$ the Bernoulli space $M^{\Z}$.
The analogous problem on the half-line $\mathbb{N}$ is also considered.
Let $A: \mathcal{B}_i \rar \R$ be an {\it observable}.
Given a temperature $T$, we analyze the main properties of the Gibbs state $\hat{\mu}_{\frac{1}{T} A}$.

In order to do our analysis we consider the Ruelle operator associated to $\frac{1}{T} A$,
and we get in this procedure the main eigenfunction $\psi_{\frac{1}{T} A}$.
Later, we analyze selection problems when the temperature goes to zero:
a) existence, or not, of the limit
$V:=\lim_{T\to 0} T\, \log(\psi_{ \frac{1}{T} A})$,
a question about selection of subactions, and,
b) existence, or not, of the limit
$\tilde{\mu}:=\lim_{T\to 0} \hat{\mu}_{\frac{1}{T}\, A}$, a question about selection of measures.

The existence of subactions and other properties of Ergodic Optimization are also considered.

The case where the potential depends just on the coordinates $(x_0,x_1)$ is carefully analyzed.
We show, in this case, and under suitable hypotheses, a Large Deviation Principle, when $T\to 0$, graph properties, etc...
Finally, we will present in detail a result due to A. C. D. van Enter and W. M. Ruszel, where the authors show,
for a particular example of potential $A$, that the selection of measure $\hat{\mu}_{\frac{1}{T}\, A} $ in this case,
does not happen.

\end{abstract}

\vspace{0.4cm}

\section{Introduction}\label{intro}

Let $(M,d)$ be a connected and compact manifold.
  We denote by $\mathcal{B}$ the Bernoulli space $M^{\N}$ of sequences represented by
$x=(x_0,x_1,x_2,x_3,....)$, where $x_i, i\geq 0$ belongs to the space (alphabet)
$M$.
By Tychonoff´s Theorem of compactness, we know   $\mathcal{B}$
is a compact metric space when equipped with the distance  given by
$d_c(x,y)=\sum_{k\geq 0} \frac{d(x_k,y_k)}{c^k}$, with $c>1$. The topologies generated by $d_{c_1}$ or $d_{c_2}$ are the same. We denote $d$ when we choose $c=2$.
In several of our results $M$ is the interval $[0,1]$ or the one-dimensional circle $\mathbb{S}^1$.

The shift $\sigma$ on ${\cal B}$ is defined by $\sigma( (x_0, x_1,x_2,x_3,....)) = (x_1,x_2,x_3,x_4,....)$. It is a continuous function on ${\cal B}$.

Let $A: \mathcal{B} \rar \R$ be an {\it observable} or {\it potential} defined on the Bernoulli space $\mathcal{B}$,
i.e. a real-valued function defined on $\mathcal{B}$. 
The potential $A$ describes an interaction between sites in the one-dimensional lattice $M^\mathbb{N}$.

For most of the results we consider here we will require $A$ to be H\"older-continuous,
which means there exist  constants $0<{\alpha}<1$ and $Hol_A>0$ such that $|A(x)-A(y)| \leq Hol_A d(x,y)^{\alpha}$. We call $\alpha$ the exponent of $A$ and $Hol_A$ the constant for $A$.  We will be interested here in the Gibbs state $\mu_A$ associated to such $A$, which will be a probability measure on ${\cal B}$. Note that the set of probability measures on ${\cal B}$ is compact for the weak* topology, 
(which is given by a metric).

 For each value $\beta=1/T$, where $T$  is the temperature, we can  consider the Gibbs state $\mu_{\beta A}$, and, we want to show in a particular  example (introduced by  A. C. D. van Enter and W. M. Ruszel \cite{van}), that there is no limit (in the weak* topology) of the family $\mu_{\beta A}$, when $\beta \to \infty$. We will present here in section \ref{sec5} all the details of the proof of this non-trivial result.

We point out that by trivial modification of the metric a Holder potential can be considered a Lipschitz potential (with no change of the topology).
Therefore, we can state our results in either case.
The assumption of $A$ being Lipschitz means that there is fast rate of decay of influence of the  potential if we are far away in the lattice.

The case of the lattice $\mathbb{Z}$, that is ${\cal B}_i=M^\mathbb{Z}$ can be treated in a similar way:
Let $A: M^\mathbb{Z}\to \mathbb{R}$ be a Lipschitz potential, and  denote by $\hat{\sigma}$ the left-shift on ${\cal B}_i$.
Any Lipschitz potential on ${\cal B}_i$ is $\hat{\sigma}$-cohomologous to a potential on ${\cal B}$
(same proof as in Proposition 1.2 \cite{PP} or, in \cite{Bow}).
We will explain this more carefully later.
To consider $\hat{\sigma}$-invariant probability measures on ${\cal B}_i$ means that the position $0\in \mathbb{Z}$
in the lattice is not distinguished (which  in general  makes sense).

We call general one-dimensional $XY$ model the setting described above. A particularly interesting case is
when we consider $M=\mathbb{S}^1$ (the unit one-dimensional circle) \cite{FH} \cite{L} \cite{van}.
This one-dimensional continuous Ising model is another important example that can be treated in the setting.
Below in section \ref{sec1} our results are for the general case of any $M$ as above.

We say that the potential $A:{\cal B} \to \mathbb{R}$ depends on the first two coordinates if $A(x)=A( x_0,x_1,x_2,..) = A( x_0,x_1)$, for any $x=( x_0,x_1,x_2,..).$ In this case $A$ is always Lipschitz. Such kind of potentials are sometimes called nearest neighbor interaction potentials.
The so-called one-dimensional $XY$ model in most of the cases assumes that $A$ depends on the first two coordinates \cite{FH}.
Special attention to this case  will be given in  section \ref{sec3}.
For example, in \cite{FH} \cite{EFS}
$$A(x)= A(x_0,x_1)= \cos(x_1-x_0-\alpha) + \gamma \cos(2\, x_0),$$
where $\alpha$ and $\gamma$ are constants.
The part $\gamma \cos(2\, x_0)$ corresponds to the magnetic term while $\cos(x_1-x_0-\alpha)$ corresponds to the interaction term.

We point out that this point of view of getting a coboundary and the systematic use of the Ruelle operator is the
Thermodynamical Formalism setting (see \cite{PP}).
This, in principle, is different from the point of view more commonly used in Statistical Mechanics on general lattices where the Gibbs measures are defined by means of a specification, DLR formalism, limit of probabilities on finite boxes (see \cite{G}, \cite{EFS}, \cite{Mayer}). We briefly address this question for a potential which depends on two coordinates in section \ref{sec-dic}.

In the Classical Thermodynamic Formalism one usually considers $M=\{1,2,...,d\}$ \cite{PP} \cite{Kel}. Here $M$ is a compact manifold with a volume form. We point out that we will use the following notation: we call a {\bf Gibbs} probability measure for $A$ the measure which is derived from a Ruelle operator, and we call the {\bf equilibrium} probability measure for $A$ the one which is derived from a maximization of Pressure (which requires one to be able to talk about entropy). We will be interested here in {\bf Gibbs} states because we need to avoid to talk about entropy. Note that the shift acting on $M^{\N}$ is such that each point has an uncountable number of pre-images. Just in some late sections  we will speak about "entropy" and  "pressure" of the potential $A$ (in general in the case it depends on two coordinates).

An interesting discussion about the several possible approaches (DLR, Thermodynamic limit in finite boxes, etc..)  to Statistical Mechanics in the one-dimensional lattice appears in \cite{Sa}.

Some of the results presented here will be used in a future related paper  \cite{LM1}.

We point out that the understanding of Statistical Mechanics   via the Ruelle Operator (Transfer Operator) allows one to get eigen-functions, and, in the limit (in the logarithm scale), when temperature goes to zero, the subaction. This helps in getting Large Deviation properties of Gibbs states when temperature goes to zero \cite{BLT} \cite{LM1}.


In the first part of this paper we describe the theory for case of general $A$ (section \ref{sec1} for positive temperature and section \ref{sec2} for zero temperature).
Later (in section \ref{sec3}) we will focus on the case where the potential $A$ depends only on the first two coordinates.
Section \ref{sec-dic} compares the setting of Thermodynamical Formalism with DLR Formalism. These two sections will
help a better understanding of Section \ref{sec5} where we present a detailed explanation of an example \cite{van} where there is
no selection of measures.

\vspace{0.6cm}

\section{Positive temperature: a generalized Ruelle-Perron-Frobenius Theorem}\label{sec1}

Let  $\mathcal{C}$ be the space of continuous functions from $\mathcal{B}=M^\mathbb{N}$ to $\R$.
We are interested in the Ruelle operator on $\mathcal{C}$  associated to the Lipschitz observable $A:M^\mathbb{N} \to \mathbb{R}$, which acts on $\psi \in \mathcal{C}$,
and sends it to $\LLW(\psi) \in \mathcal{C}$ defined by
$$\LLW(\psi)(x)=\II e^{A(ax)}\psi(ax)\, d\,a\,,
$$
for any  $x=(x_0,x_1,x_2,....)\in \mathcal{B}$,
where $ax$ represents the sequence
$(a,x_0,x_1,x_2,....)\in \mathcal{B}$, and $d\,a$ is the Lebesgue probability measure on $M$. Note that $\sigma(ax)=x$.

A major difference between the settings of the Classical Bowen-Ruelle-Sinai Thermodynamic Formalism setting
and the $XY$ model is that here, in order to define the Ruelle operator, we need an a priori measure
(for which we consider in most of the cases the  Lebesgue probability measure $da$ on $\mathbb{S}^1$).

Some of the results of the present section are generalization of theorems in \cite{LMST}.

The operator $\LLW$ will help us to find the Gibbs state for $A$. First we will show the existence of a main eigenfunction for ${\cal L}_A$, when $A$ is Lipschitz. Part of our proof follows  the reasoning of section 7 in \cite{Bousch-walters} (which considers $M=\{1,2,..,d\}$), adapted to the present case.

We begin by defining another operator on $\mathcal{C}$. Let $0<s<1$, and define, for $u \in \mathcal{C}$, $\TT(u)$ given by
$$\TT(u)(x) = \log\left( \II e^{A(ax)+su(ax)}\, da \right).$$

\begin{proposition}
If $0<s<1$ then
$\TT$ is an uniform contraction map.
\end{proposition}

{\it Proof.:}
$$\left|\TT(u_1)(x)-\TT(u_2)(x)\right|=
\left|\log\left( \frac{ \II e^{A(ax)+su_1(ax)}}{ \II e^{A(ax)+su_2(ax)}}
\right)\right|=$$
$$=\left|\log\left( \frac{ \II e^{A(ax)+su_2(ax)+ su_1(ax)  -su_2(ax)}
}{ \II e^{A(ax)+su_2(ax)}}
\right) \right| \leq$$
$$\leq \log\left( \frac{ \II e^{A(ax)+su_2(ax)+
s\|u_1  -u_2\|
}
}{ \II e^{A(ax)+su_2(ax)}}
\right) = s \|u_1  -u_2\| \,.$$

\cqd

\vspace{0.3cm}

Let $u_s$ be the unique fixed point for $\TT$. We have
\begin{equation}\label{logeigforS}
\log\left(\II e^{A(ax)+su_s(ax)}\, da\right)=u_s(x)\,.
\end{equation}

\begin{proposition}\label{famequic}
The family $\{ u_s\}_{0<s<1}$ is an equicontinuous family of functions.
\end{proposition}

{\it Proof.:} 
 Let $H_s(x,y)=u_s(x)-u_s(y)$. By \eqref{logeigforS} we have
\begin{eqnarray*}e^{u_s(x)}&=&\II e^{A(ax)+su_s(ax)}\\
&=&\II e^{A(ay)+su_s(ay)}e^{A(ax)-A(ay)+s[u_s(ax)- u_s(ay)]}\\
&\leq& e^{u_s(y)}\max_{a}\{{e^{A(ax)-A(ay)+s[u_s(ax)- u_s(ay)]}}\}.
\end{eqnarray*}
Hence

$$e^{u_s(x)-u_s(y)} \leq  \max_{a}\{{e^{A(ax)-A(ay)+s[u_s(ax)- u_s(ay)]}}\}, $$ and this implies

$$H_s(x,y)=u_s(x)-u_s(y) \leq  \max_{a}[A(ax)-A(ay)+sH_s(ax,ay)]. $$

Proceeding by induction we get
$$H_s(x,y)\leq \max_{\theta\in \mathcal B}\sum_{n=0}^{\infty}s^n[A(\theta_n....\theta_0x)-A(\theta_n...\theta_0y)   ]\leq  $$
$$  \leq Hol_A \max_{\theta\in \mathcal B}\sum_{n=0}^{\infty}s^n d((\theta_n....\theta_0x),(\theta_n...\theta_0y))^{\alpha}   \leq  $$
 $$ \leq Hol_A \sum_{n=0}^{\infty}\left(\frac{s}{2^{\alpha}}\right)^nd(x,y)^{\alpha}\leq
 \frac{2^{\alpha}}{2^{\alpha}-1} Hol_A d(x,y)^{\alpha}\,.$$

{\bf Remark 1:}
This shows that $u_s$ is Lipschitz, and, moreover, that $u_s$, $0\leq s<1$, is an equicontinuous family.
Note the very important point: the Lipschitz constant of $u_s$, is given by  $\frac{2^{\alpha}}{2^{\alpha}-1} Hol_A$, and depends only on the Holder constant for $A$,
but does not depend on $s$.

\cqd

\bigskip


Let
$$S_n(z)=S_{n,A}(z)=\sum_{k=0}^{n-1} A\circ \sigma^k(z)\,.$$

Note that iterates of the operator ${\cal L}_{ A}$ can be written with the use of $S_{n,A}(z)$.
$$ {\cal L}_{ A}^n (w)(x)= \int_{{\mathbf a} \in M^n}\, e^{S_{n,  A}\,({\mathbf a}x)} w({\mathbf a}x) \,d{\mathbf a} .$$


\medskip

\begin{theorem}\label{RPF_eig_nonnorm}
There exists a strictly positive Lipschitz eigenfunction
 $\evc$
for $\LLW: {\cal C} \to {\cal C}$ associated to a strictly positive eigenvalue
$\evl$. The eigenvalue is simple and it is equal to the spectral radius.
\end{theorem}

\noindent{\bf Proof.}  It follows from the fixed point equation that for any $x$
 $$-||A||+ s \min u_s \leq u_s(x)\leq ||A||+ s \max u_s.$$

 Therefore, $-||A|| \leq (1-s) \min u_s\leq (1-s) \max u_s\leq||A||$, for any $s$. Consider a subsequence $s_n \to 1$ such that $ [\,(1-s_n) \, \max u_{s_n}\,]\,\to k$.

 The family $\{ u^*_s=u_s-\max u_s \}_{0<s<1}$ is equicontinuous and uniformly bounded.

 Therefore, by Arzela-Ascoli $\{u_{s_n}^*\}_{n\geq 1}$  has an accumulation point in $\mathcal{C}$, which we will call $u$.

Observe that for any $s$
$$e^{u^*_s (x)}= e^{u_s(x) - \max u_s}= $$
$$e^{-(1-s) \max u_s + u_s (x) - s \max u_s}=$$
$$e^{-(1-s) \max u_s}\,\int e^{  A(ax)+ (s  u_s (ax) - s \max u_s)}\, da.$$

Taking limit where $n$ goes to infinity for the sequence $s_n$ we get that $u$ satisfies
$$e^{u(x)} =e^{-k} \,  \int e^{  A(ax)+   u (ax)}\,da.$$

In this way we get a {\bf positive} Lipschitz eigenfunction $\psi_{A}=e^{u}$ for ${\cal L}_A$ associated to the eigenvalue $\lambda_A=e^{k}$.
\cqd


{\bf Remark 2:}
To prove that $u$ is Lipschitz, we just
use the fact that $u$ is the limit of a sequence of uniformly Lipschitz functions (i.e. Lipschitz functions with same Lipschitz constant).  Using that $u$ is a bounded function we have that
$\psi_A=e^u$ is also Lipschitz. Note a very important point: the Lipschitz constant of $u=\log(\psi_A)$
is given by $\frac{2^{\alpha}}{2^{\alpha}-1} Hol_A$
 (see Remark 1  in the end of the proof of  Proposition \ref{famequic}).

The property that the eigenvalue is simple and maximal  follows from the same reasoning as in page 23 and 24 of \cite{PP}.
For example, to prove that the eigenvalue is simple we suppose there are two eigenfunctions $\psi_1$ and $\psi_2$. Let $t=\min\{\psi_1 / \psi_2\}$. Then $\psi_3=\psi_1-t\psi_2$ is a non-negative eigenfunction which vanishes at some point $z\in \mathcal{B}$. Therefore
$$0=\lambda_A^n \psi_3(z)
=\int_{{\mathbf a} \in M^n}\, e^{S_{n,  A}\,({\mathbf a}z)} \psi_3({\mathbf a}z) \,d{\mathbf a}\,,$$
which implies $\psi_3({\mathbf a}z)=0 \,\,\forall \,{\mathbf a} \in M^n$, $\forall n$, which makes $\psi_3=0$.

\cqd

\vspace{0.3cm}

Note that
\begin{equation}\label{normalizado}\II \frac{e^{A(ax)}\evc(ax)}{\evl \evc(x)}da=1\,,\,\forall x \in \mathcal{B}\,.\end{equation}
If a potential $B$ satisfies
$$\II e^{B(ax)}da = 1 \,,\,\forall x \in \mathcal{B}\,,$$
 which means $\mathcal{L}_B(1)=1$, we say that
$B$ is normalized.

Let $$\bar A = A+\log  \evc-\log \evc \circ \sigma -\log \evl,$$
where $\sigma:\BB \rar \BB$ is the usual shift map.
Equation \eqref{normalizado} shows that
$\bar A$ is normalized. It is also Lipschitz (Holder). In this case the main eigenvalue is $1$ and the main eigenfunction is constant equal  to $1$
(in fact we can prove, using proposition \ref{RPF_normalized}, that there is only one strictly positive eigenfunction, the one associated to the maximal eigenvalue).


Remember that, given $x=(x_0,x_1,x_2,...)\in {\cal B}$ and $a\in M$, we denote by $ax\in {\cal B}$ the element  $ax=(a,x_0,x_1,x_2,...),$ i.e., any $y\in {\cal B}$ such that $\sigma(y)=x$ is of this form.

We define the Borel sigma-algebra ${\cal F}$ over $\BB$ as  the $\sigma$-algebra
generated by the {\it cylinders}. By this we mean the sigma-algebra generated by sets of the form
$B_1\times B_2\times\,...\,\times B_n \times M^{\mathbb{N}}$, where $n \in \mathbb{N}$, and $B_j, j \in \{1,2,..,n\}$, are open sets in $M$.
Similar definitions can be considered for ${\cal B}_i$.

We say a probability measure $\mu$ over ${\cal F}$ is invariant, if for any Borel set $B$, we have that $\mu(B)= \mu(\sigma^{-1} (B)).$ This corresponds to stationary probability measures for the underlying stochastic process $X_n$, $n \in \mathbb{N}$, with state space $M$. We denote by ${\cal M}_\sigma$ the set of invariant probability measures. Similar definitions can be considered for ${\cal B}_i$.

We present below a generalization of results considered in \cite{PP}.

We  define the dual operator $\DLLW$ on the space of the Borel measures on $\BB$ as the operator that sends a measure $v$ to the measure $\DLLW(v)$  defined by
$$\int_{\BB} \psi\, d\DLLW(v) =  \int_{\BB} \LLW(\psi)\, dv\,.
$$
for any $\psi \in \mathcal{C}$.

Now we want to find an eigen-probability for $\LLW^*$. This will help us to find the Gibbs state for the potential $A$.

\begin{proposition}\label{RPF_normalized}
If the observable $\bar A$ is normalized, then there exists an unique fixed point $m=m_{\bar A}$ for $ {\cal L }_{\bar A}^*$. Such a probability measure $m$ is $\sigma$-invariant, and for all
Holder continuous function $\omega $
we have that, in the uniform convergence topology,
$${\cal L }_{\bar A}^n\omega \rar \int_{\BB} \omega dm\,.$$
\end{proposition}
Here ${\cal L }_{\bar A}^n$ denotes the $n$-th iterate of the operator ${\cal L }_{\bar A}:{\cal C}\to {\cal C}$.

\medskip

{\it Proof.:} We begin by proving that the normalization property implies that the convex and compact set of Borel probability measures on $\BB$ is preserved by the operator ${\cal L}_{\bar A}^*$: in order to see that, note that for $\mu$ a Borel probability measure on $\BB$, we have
$${\cal L}_{\bar A}^*(\mu)(\BB) = \int_{\BB} 1 \,\,d{\cal L}_{\bar A}^*(\mu) =
 \int_{\BB}  {\cal L}_{\bar A} (1) d\mu  =
 \int_{\BB}  1 \,\,d\mu= \mu(\BB) =1$$ where the 
 third equality is precisely the normalization hypothesis.

By the Tychonoff-Schauder theorem let $m$ be a fixed point for the operator ${\cal L}_{\bar A}^*$.

To prove that $m$ is $\sigma$-invariant, we begin by observing that
$${\cal L}_{\bar A}(\psi \circ \sigma)(x)=
\II e^{\bar A (ax)}\psi\circ \sigma(ax)da = \II e^{\bar A (ax)}\psi(x)da = \psi(x).$$
Note that the normalization hypothesis is used in the last equality.

Therefore,  if $\psi \in \mathcal{C}$, then
$$\int_{\BB} \psi \circ \sigma dm =\int_{\BB} \psi \circ \sigma d{\cal L}_{\bar A}^*(m) = \int_{\BB} {\cal L}_{\bar A}( \psi \circ \sigma) dm =\int_{\BB}  \psi dm .
$$
which implies the invariance property of $m$.

\bigskip


Before finishing the proof of proposition \ref{RPF_normalized}, we will need two claims. The first is a special estimate which will be important in the rest of this section.

\medskip
{\it Claim:} For any  Holder potential $A$, if $\|w\|$ denotes the uniform
norm of the Holder function $w:{\cal B}\to \mathbb{R}$,
 we have 
$$|\LLW^n (w) (x) - \LLW^n(w) (y) |\leq \left[C_{e^A}\|w\| \left( \frac{1}{2^{\a}}+...+\frac{1}{2^{n\a}}\right)+\frac{C_w}{2^{n\a}} \right]
d(x,y)^{\a},$$
where $C_{e^A}$ is the Holder constant of $e^A$ and $C_w$ is the Holder constant of $w$.

{\it Proof of the Claim:} :
We prove the claim by induction. Suppose $n=1$. We have

$$| \LLW(w)(x)-\LLW(w)(y)| \leq
$$
$$ \leq \int_{M} |e^{A(ax)}-e^{A(ay)}| \cdot |w(ax)| da + \int_{M}e^{A(ay)} |w(ax)-w(ay)|da \leq $$
$$ \leq (C_{e^A} \|w\| + C_w )\frac{d(x,y)^{\a}}{2^{\a}},  $$
where in the last inequality we used the normalization property of $A$.
In particular we can say that the Holder constant of $\LLW(w)$ is given by
\begin{equation}\label{holdconstantLw}
C_{\LLW(w)} = \frac{C_{e^A} \|w\| + C_w  }{2^{\a}}\,.
\end{equation}

Now, suppose the Claim holds for $n$. We have
$$|\LLW^{n+1} (w) (x) - \LLW^{n+1}(w) (y) | =
|\LLW^n (\LLW(w)) (x) - \LLW^n(\LLW(w)) (y) |\leq
$$
$$ \leq \left[C_{e^A}\|\LLW(w)\| \left( \frac{1}{2^{\a}}+...+\frac{1}{2^{n\a}}\right)+\frac{C_{\LLW(w)}}{2^{n\a}} \right]
d(x,y)^{\a},
$$
and, therefore the claim is proved when we use \eqref{holdconstantLw} and  $\|\LLW(w)\|\leq\|w\|$ which is consequence of the normalization property of $A$.
\bigskip

As a consequence, the set $\{{\cal  L}_{\bar A}^n \omega\}_{n\geq 0}$ is equicontinuous.
In order to prove that $\{{\cal  L}_{\bar A}^n \omega\}_{n\geq 0}$ is uniformly bounded we  use again the normalization condition which implies $\|{\cal  L}_{\bar A}^n \omega\|\leq\|w\| \,,\forall n \geq 1$.


\bigskip
By the Arzela-Ascoli Theorem let $\bar \omega$ be an accumulation point for  $\{{\cal L}_{\bar A}^n \omega\}_{n\geq 0}$, i.e., suppose there exists a subsequence $\{n_k\}_{k\geq 0}$ such that $$\bar \omega (x) = \lim_{k \geq 0} {\cal L}_{\bar A}^{n_k}\omega(x)\,.$$

{\it Second Claim:} : $\bar \omega$ is a constant function.

 The proof of this second claim is similar to the reasoning of page 25 \cite{PP}.

Now that $\bar \omega$ is a constant function we can prove  that
$$ \bar \omega = \int_{\BB} \bar \omega dm = \lim_k \int_{\BB} {\cal L}_{\bar A}^{n_k} \omega dm
= \lim_k \int_{\BB}  \omega d ({\cal L}_{\bar A}^*)^{n_k}(m)
= \int_{\BB}  \omega dm,$$
which shows that $\bar \omega$ does not depend on the subsequence chosen. Therefore, for any $x \in \BB$ we have $${\cal L}_{\bar A}^{n} \omega (x) \rar \bar \omega = \int_{\BB} \omega dm\,.$$
The last limit shows that the fixed point $m$ is unique.

\cqd

\vspace{0.3cm}

\begin{proposition}\label{eigenprob}
Let $A$ be a Holder, not necessarily normalized potential, and $\psi_A$ and $\lambda_A$ the eigenfunction and eigenvalue given by  theorem \ref{RPF_eig_nonnorm}.
To the potential $A$ we associate the normalized potential
$\bar A = A+\log  \evc-\log \evc \circ \sigma -\log \evl$. Let $m$ be the unique probability measure that satisfies ${\cal L}_{\bar A}^*(m)=m$, given by proposition
\ref{RPF_normalized}.

(a) the measure $$\rho_A= \frac{1}{\psi_A} \,\, m$$ satisfies ${\cal  L}_A^* (\rho_A)= \lambda_A \rho_A$. Therefore, $\,\rho_A$  is an eigen-probability for ${\cal  L}_A^*$.

(b) for any Holder  $\phi:{\cal B} \to \mathbb{R}$, we have that
 $$\frac{{\cal L}_A^n (\phi)}{\lambda_A^n} \to \, \psi_A \int\, \phi \,d \rho_A.$$
\end{proposition}

\medskip
{\it Proof:} (a)
${\cal L}_{\bar A}^*(m)=m$ implies that for any $\psi \in \mathcal{C}$, we have
\begin{eqnarray*}
  \int \psi dm &=& \int \psi d {\cal  L}_{\bar A}^*(m)\\
   &=& \int  {\cal  L}_{\bar A}(\psi) dm \\
   &=& \int \left( \int \psi(ax) e^{\bar A(ax)}da  \right) dm(x) \\
   &=& \int \left( \int \psi(ax) \frac{e^{A(ax)}  \psi_A(ax)}{\lambda_A \psi_A(x)} da  \right) dm(x) \,.\\
\end{eqnarray*}
Now, if $\varphi \in \mathcal{C}$, making $\psi=\frac{\varphi}{\psi_A}$ in the last equation we have
$$\int \frac{\varphi}{\psi_A} dm = \frac{1}{\lambda_A} \int \left( \int \varphi(ax) \frac{e^{A(ax)} }{ \psi_A(x)} da  \right) dm(x)\,,$$
which is equivalent to
\begin{equation}\label{gibbs_eig}
    \lambda_A \int \varphi d \rho_A = \int {\cal  L}_A(\varphi)  d \rho_A \,
\end{equation}
or $${\cal  L}_A^* (\rho_A)= \lambda_A \rho_A\,.$$

(b) We have that $A = \bar A - \log  \evc + \log \evc \circ \sigma +\log \evl,$ and therefore

$$S_{n,A}(z) \equiv \sum_{k=0}^{n-1} A\circ \sigma^k(z)=
S_{n,\bar A}(z) -\log \psi_A + \log \psi_A \circ \sigma^n + n \log \lambda_A
\,,$$
which makes 
$$\frac{{\cal L}_A^n(\phi)(x)}{\lambda_A^n}=
\frac{1}{\lambda_A^n}\int_{{\mathbf a}\in M^n} e^{ S_{n,A}({\mathbf a}x)} \phi({\mathbf a}x) d{\mathbf a}=$$
$$=\psi_A(x) \int_{{\mathbf a}\in M^n} \frac{e^{ S_{n,\bar A}({\mathbf a}x)}}{\psi_A({\mathbf a}x)} \phi({\mathbf a}x) d{\mathbf a}=$$
$$=\psi_A(x) {\cal L}_{\bar A}^n\left(\frac{\phi}{\psi_A}\right)
\to \psi_A(x) \int \frac{\phi}{\psi_A} dm_{\bar A}$$
where the convergence on $n$ in the last line comes from
Proposition \ref{RPF_normalized}. 
\cqd

\medskip
{\bf Remark 3:} From now on we will call
$m_{\bar A}$ the eigen-probability for ${\cal  L}_{\bar A}^*$.
One can show that the eigen-probability
$\rho_A = \frac{1}{\psi_A} \,\, m_{\bar A}$ is the unique eigen-probability for ${\cal  L}_A^*$.
Also,  it is not necessarily invariant for the shift $\sigma$.


{\bf We call $m_{\bar A}$ the Gibbs state for $A$.} 
This probability measure $m_{\bar A}$ over ${\cal B}$ is invariant for the shift and describes the statistics of the interaction described by $A$.
It is usual to call the probability measure $m_{\bar A}$ the Gibbs state (in the Thermodynamic Formalism setting \cite{PP}) for the interaction given by $A$.


We point out that the probability measure $\rho_A$ is positive on open sets of ${\cal B}$. Suppose the metric space $M=\mathbb{S}^1$. The projection of this probability measure on the first two coordinates $\mathbb{S}^1\times \mathbb{S}^1$ is absolutely continuous with respect to Lebesgue probability measure on $\mathbb{S}^1\times \mathbb{S}^1$.
This is so because, if $B$ is Borel in $[0,1]^2$, then from \eqref{gibbs_eig} we have
$$\int I_{(x_0, x_1)\in B} \, d \rho_A= \frac{1}{\lambda_A^2}\int  {\cal L}_{\bar A}^2\,(I_{(x_0, x_1)\in B}) \, d \rho_A,$$
and, for any $x\in{\cal B}$
$$ {\cal L}_{\bar A}^2 ( I_{(x_0, x_1)\in B})(x)= \int_M \int_M\, e^{S_{2, \bar A}\,(abx)} I_{(x_0, x_1)\in B}(abx) \,da \, db.$$

{\bf Remark 4:}
 If we consider instead a Holder potential $B:{\cal B}_i=M^{\mathbb{Z}} \to \mathbb{R}$, where
 $${\cal B}_i=\{( ...,x_{-2}, x_{-1}, x_0, x_1,x_2,...)|\, x_i \in M, i \in \mathbb{Z}\},$$
 then, we first derive  (as in Proposition 1.2 \cite{PP} or, in \cite{Bow}) the associated cohomologous Holder potential $A: {\cal B} \to \mathbb{R}$ (the Holder class can change), then proceed as above to get $\rho_A$ over ${\cal B}$. Finally, we consider the natural extension $\hat{\rho}_A$ of $\rho_A$ on ${\cal B}_i$
 (see \cite{PY} \cite{Bow}),
 and we solve in this way the Statistical Mechanics problem for the interaction described by $B$ in the lattice $\mathbb{Z}$: it´s the probability measure  $\hat{\rho}_{\beta A}.$

 Note that if $C$ is a set that depends just on the coordinates $x_0,x_1$, then $\rho_{\beta A}(C)= \hat{\rho}_{\beta A}(C).$ For sets $C \subset {\cal B}_i$, of this form, we can use  without loss of generality  $\rho_{\beta A}(C)$ or   $\hat{\rho}_{\beta A}(C).$

\begin{proposition}\label{proponly}

The only Lipschitz continuous eigenfunction $\psi$ of ${\cal L}_A$ which is totally positive is $\psi_A$ (the one associated to the maximal eigenvalue $\lambda_A$).

\end{proposition}

{\it Proof:} Suppose $\psi:{\cal B} \to \mathbb{R}$ is a Lipschitz continuous  eigenfunction  of ${\cal L}_A$ associated to the eigenvalue $\beta$.

It follows from the above that
$\frac{{\cal L}_A^n (\psi)}{\lambda_A^n} \to \, \psi_A \int \psi d \rho_A$, when $ n \to \infty$.

 Therefore, if $\psi >c>0$, then $\int \psi d \rho_A>0$. Moreover, ${\cal L}_A^n (\psi) = \beta^ n \psi$. This is only possible if $\beta=\lambda_A$ and $\psi=\psi_A$.

\cqd

\vspace{0.3cm}

It is easy to see that if $A$ is Holder with exponent $\alpha$, and, denoting
${\cal H}_\alpha$, the set of real-valued functions with Holder exponent $\alpha$, then ${\cal L}_{\bar A} : {\cal H}_\alpha\to {\cal H}_\alpha.$

For $w \in {\cal H}_\alpha$, denote $|w|_\alpha= \sup_{x \neq y} \frac{|w(x)- w(y)|}{d(x,y)^\alpha}$. It is known that ${\cal H}_\alpha$ is a Banach space for the norm
$$ ||w||_\alpha= | w|_\alpha + ||w||,$$
where $||w||$ is the uniform norm of $w$.

When $\alpha=1$ we are considering the space of Lipschitz functions ${\cal H}_1$.

We note that $\mathcal{K}_{\alpha}\equiv\{w \in {\cal H}_\alpha\,,\, ||w||_\alpha\leq 1\}$ is compact in the uniform norm as a subset of ${\cal C}.$ To prove that, we just need to observe that the definition of the norm $||w||_\alpha$ implies that $\mathcal{K}_{\alpha}$ is a equicontinuous and uniformly bounded set, and then we have the result directly by using Arzela-Ascoli´s theorem.

We can also prove that $\mathcal{K}_{\alpha}^{ A}\equiv\{w \in {\cal H}_\alpha\,,\, \int_w dm_{ A}=0\,\,,\,\,||w||_\alpha\leq 1\}$ is compact in the uniform norm. 
 For doing that, let
$I_{m_{ A}}: {\cal H}_\alpha \to \mathbb{R}$ be given by $I_{m_{ A}}(w)=\int w dm_A$. We have that
$I_{m_{ A}}$ is a bounded linear operator, and therefore $I_{m_{ A}}^{-1} \{0\}$ is a closed  subset of ${\cal H}_{\alpha}$. Now $\mathcal{K}_{\alpha}^{ A}= \mathcal{K}_{\alpha} \cap I_{m_{A}}^{-1} \{0\}$ is compact. 

\begin{proposition}\label{espectro} Suppose $\bar A$ is normalized, then
the eigenvalue $\lambda_{\bar A}=1$ is maximal. Moreover, the remainder of the spectrum of ${\cal L}_{\bar A}: {\cal H}_\alpha\to {\cal H}_\alpha$ is contained in a disk centered at zero with radius strictly smaller than one.
\end{proposition}

\noindent{\bf Proof.}
Remember that $1$ is the eigenfunction associated to the eigenvalue $1$.
We will show that ${\cal L}_{\bar A}$ restricted to
$\mathcal{K}_{\alpha}^{\bar A}$ has spectral radius strictly smaller than $1$.
We know from proposition \ref{RPF_normalized} that ${\cal L}^{k}_{\bar A}$ converges to zero in the compact set $\mathcal{K}_{\alpha}^{\bar A}$.

The normalization hypothesis implies $||\mathcal{L}_{\bar A}^{n+1}(w)|| \leq ||\mathcal{L}_{\bar A}^{n}(w)||\,\forall n \geq 0$.
We will now prove that this monotonicity property implies that the convergence above is uniform. 
More precisely, we have

{\it Claim:}
Given a small $\epsilon$, there exists $N=N_{\epsilon}\in \mathbb{N}$ such that

$$||\mathcal{L}_{\bar A}^{n}(w)|| < \epsilon \,\,\forall \, n\,\geq N\,,\,\forall w \in \mathcal{K}_{\alpha}^{\bar A}\,.$$

To prove this claim, let $C_n\equiv \{w \in \mathcal{K}_{\alpha}^{\bar A} \,:\, ||\mathcal{L}_{\bar A}^{m}(w)|| < \epsilon \,\forall\, m \geq n   \}$.
The monotonicity property implies $C_n \subseteq C_{n+1}$ and also that $C_n$ is an open set in the uniform norm, while ${\cal L}^{k}_{\bar A}(w)\to 0$ implies $\cup_n C_n = \mathcal{K}_{\alpha}^{\bar A}$. Therefore, compactness of $\mathcal{K}_{\alpha}^{\bar A}$ implies $\mathcal{K}_{\alpha}^{\bar A}=C_N$ for some $N\in \mathbb{N}$.

The last claim is easy to prove and can be enunciated as:

{\it Claim:} There exists $C>0$ such that

$\forall n\in \mathbb{N}$ and $w \in \mathcal{H}_{\alpha}$
$$|{\cal L}^n_{\bar A}(w)|_\alpha\leq C ||w|| + \frac{|w|_\alpha}{(2^\alpha)^n}.$$

Now, for any given $n$ and $k$, using the last Claim  we have for $w \in \mathcal{H}_{\alpha}$
$$|{\cal L}^{n+k}_{\bar A}(w)|_\alpha\leq C ||{\cal L}^{k}_{\bar A}(w)|| + \frac{|{\cal L}^{k}_{\bar A}(w)|_\alpha}{(2^\alpha)^n}\leq C||{\cal L}^{k}_{\bar A}(w)|| + C\, \frac{||w||}{(2^\alpha)^n} +\, \frac{|w|_\alpha}{(2^\alpha)^{n+k}}.$$

Therefore, if $\epsilon$ is small enough and $n\geq N_{\epsilon}$, we have that for all $w\in \mathcal{K}_{\alpha}^{\bar A}$
$$||{\cal L}^{n+k}_{\bar A} (w)||_\alpha \leq \epsilon<1.$$

In this case the spectral radius is smaller than $\epsilon^{\frac{1}{n+k}}.$

\cqd

\vspace{0.3cm}

We denote $\lambda^1_{\bar A}< \lambda_{\bar A}$
the spectral radius of ${\cal L}_{\bar A}$ when restricted to the set  $\{w\in {\cal H}_\alpha: \int w \,d m_{\bar A}=0\}$.

Now we will show the exponential decay of correlation for Holder functions.

\begin{proposition}\label{dual e correlacao}
If $v,w \in {\cal L}^2 (m_{\bar A})$ are such that $w$ is Holder and $\int w \,  d m_{\bar A}=0$, then, there exists $C>0$ such that for all $n$
$$\int (v \circ \sigma^n) \, w \, d m_{\bar A}\leq C \, (\lambda^1_{\bar A})^n$$
\end{proposition}

\noindent{\bf Proof.} This follows from
$$ \int (v \circ \sigma^n) \, w \, d m_{\bar A}= \int v  \,{\cal L}_{\bar A}^n w \, d m_{\bar A}.$$

\cqd

\vspace{0.3cm}

 The above proposition implies that $m_{\bar A}$ is mixing (same reasoning as in section 2 of \cite{Kel} which considers the case of the shift on $\{1,2,..,d\}^\mathbb{N}$).

\begin{proposition}\label{erg}

The invariant probability measure $m_{\bar A}$ is ergodic.

\end{proposition}

\noindent{\bf Proof.} If a dynamical system is mixing then it is ergodic (see section 2 in \cite{Kel}).

\cqd

\vspace{0.3cm}

A major difference of the general $XY$ Model to the Thermodynamic Formalism setting (in the sense of \cite{PP} \cite{Kel}) in $\{1,2,...,d\}^\mathbb{Z}$ is that here we can not define in the traditional way (via dynamic partitions) the concept of entropy of an invariant probability measure $\mu$ (defined over the sigma algebra ${\cal F}$ of ${\cal B}$). Each element $x\in {\cal B}$  has an uncountable set of pre-images and this is a problem.

Note that there exist invariant probabilities (for instance, singular with respect to Lebesgue measure) for the shift on ${\cal B}$ which have Kolmogorov entropy arbitrarily  large.

For the other hand, in the DLR-Gibbs theory, see  \cite{Israel}, a definition of entropy is presented and the variational principle at positive temperatures is worked out and proved.
But here we take another path, just in terms of transfer operators,  and we present the theory of the Ruelle operator for continuous-spin models, and also including some
noncontinuous potentials, which is not part of standard treatments.

Note that the Gibbs state formalism via boundary conditions, as in \cite{G}, does not require, in principle, to talk about entropy (see also our Section \ref{sec-dic}). We will address the question about entropy when the potential depends on two coordinates in section \ref{sec3}.

In Statistical Mechanics, for a fixed interaction $A$ under a certain temperature $T>0$, up to a multiplicative constant, the natural potential to be considered  is $\frac{1}{T}\,A$. We denote $\beta=\frac{1}{T}$, and, using the results above we  can consider the corresponding eigenfunction $\psi_{\beta A}$, eigenvalue $ \lambda_{\beta A}=\lambda_\beta$, and the Gibbs state which now will be denoted $\mu_{\beta A}$.

What happen with these two objects when $T\to 0$ (or, $\beta \to \infty$),  is the purpose of the next section.

\vspace{0.6cm}

\section{Zero temperature: calibrated subactions, maximizing probability measures and selection of probability measures}\label{sec2}

In this section and also in the next two sections we will consider, among other issues, questions involving selections of probability measures when the temperature goes to zero, maximizing probability measures for a given potential and existence of calibrated subactions. Among other results we will show that, under some conditions, the sequence $\{\mu_{\beta A}\}$ of Gibbs states for the potential $\beta A$ converges to a measure $\mu_{\infty}$ which has the property of maximizing the integral $\int A d\mu$  among all invariant measures $\mu$ for the shift map.
Sometimes such convergence will not occur (this is what we call non selection of probability measures - a very interesting example
due to A. C. D. van Enter and W. M. Ruszel will be presented in section \ref{sec5}).

We will also consider calibrated subactions, which is an important tool that allows one to identify the support of the maximizing probability measure $\mu_{\infty}$ (see equation \eqref{funcaoR} below), and can be used to relate the maximal eigenvalues of the Ruelle operator to the value $m(A)= \int A d\mu_{\infty}$
(see theorem \ref{existsubac}). Existence of calibrated subactions  are also related to the existence of  large deviation principles for the convergence of $\{\mu_{\beta A}\}$ to  $\mu_{\infty}$ (see theorem \ref{teoLDP} in section \ref{sec3}).

Some of the problems discussed here are usually called ergodic optimization problems (see \cite{J1}). We refer the reader to \cite{CLO} for question related to Ergodic Transport Theory.

\bigskip

Consider a fixed Holder potential $A$ and a real variable $\beta>0$.
We denote by $\psi_{\beta A}$ the eigenfunction for the Ruelle operator associated to $\beta A$.

{\bf Remark 5:}
Given $\beta$ and $A$, the Lipschitz constant of $u_\beta$, such that $\psi_{\beta A}= e^{u_\beta}$, depends on the Holder constant for $\beta\, A$ (see Remarks 1 and 2). More precisely, the Lipschitz constant of $u_{\beta}=\log(\psi_{\beta A})$
 is given by $\beta \frac{2^{\alpha}}{2^{\alpha}-1} Hol_{ A}$.
Therefore, $\frac{1}{\beta}\log(\psi_{\beta A})$, $\beta>0$, is equicontinuous. Note that it is also uniformly bounded from the reasons described below.

A possible renormalization condition for $\psi_{\beta A}$ \cite{CLT} is $\int \psi_{\beta A} \, d \rho_{\beta A}=1$, where $\rho_{\beta A}$ is the eigen-probability for ${\cal L}^*_{\beta A}$ (see proposition \ref{eigenprob} and remark 3). 
For each $\beta>0$ the normalization hypothesis $\int \psi_{\beta A} \, d \rho_{\beta A}=1$ implies the existence of  $x_\beta \in {\cal B}$ such that  $\psi_{\beta}(x_\beta)=1$. Here we are using the connectedness hypothesis of ${\cal B}$. When $\beta \to \infty$ we have that $x_{\beta_k}\to \bar x$, for a subsequence. Note that when we normalize $\psi_{\beta A}$ the Holder  constant of $\log (\psi_{\beta A})$ remains unchanged, which assures the uniformly continuous property of the family $1/\beta \log (\psi_{\beta A})\,,\beta>0$. 
 Moreover, the normalization hypothesis and Remark 5 implies that $1/\beta \log (\psi_{\beta A})\,,\beta>0$ is uniformly bounded.

Therefore, there exists a subsequence $\beta_n\to \infty$, and $V$ Lipschitz, such that on the uniform convergence
$$V:=\lim_{n\to\infty}
\frac{1}{\beta_n}\log(\psi_{\beta_n A}).$$

Consider  point $p_0\in {\cal B}$. Another possible normalization for the eigenfunction $\psi_{\beta A}$ is to assume that $\psi_{\beta A}(p_0)=1$.
We will prefer this latter form.

 By selection of a function $V$, when the temperature goes to zero (or, $\beta\to \infty$), we mean the existence of the limit (in the uniform norm)
$$V:=\lim_{\beta\to\infty}
\frac{1}{\beta}\log(\psi_{\beta A}).$$

The existence of the limit when $\beta\to \infty$ (not just of a subsequence), in the general case, is not an easy question.

In this section we denote  $\mu_{\beta A}$ the Gibbs state for the potential $\beta A$, i.e. the eigen-probability of $\mathcal{L}_{\bar A}^*$, where
$\bar A = A+\log  \evc-\log \evc \circ \sigma -\log \evl$.

By selection of a measure $\tilde{\mu}_\infty$, when the temperature goes to zero (or, $\beta\to \infty$), we mean the existence of the limit (in the weak$^*$ sense)
$$\tilde{\mu}_\infty:=\lim_{\beta\to\infty} \mu_{\beta A}.$$

In some sense $V$ is what one can get in the limit, in the $\log$-scale,  from the eigenfunction (at non-zero temperature), and $\tilde{\mu}_\infty$ is the Gibbs state at temperature zero.

Even if $A$ is Lipschitz not always the above limit on $\mu_{\beta A}$, $\beta \to \infty$, exist. In fact we will show an interesting example in  section \ref{sec5} (due to A. C. D. van Enter and W. M. Ruszel) where there is no limit for $\mu_{\beta A}$, as $\beta \to \infty$.

Some theorems in this section are generalizations of corresponding ones in \cite{LMST} (which consider only potentials $A$ which depend on two coordinates). Related results appear in \cite{GT1} and \cite{GT2}. Results about selection (or, non selection) in the setting of Thermodynamic Formalism appear in \cite{BLM}  \cite{BLL} \cite{chazottes-cv} \cite{leplaideur-max} \cite{Bremont} \cite{LM1}.

Some of the proofs and results presented in the present section are similar to other ones in Ergodic Optimization \cite{J1} and  Thermodynamic Formalism,
but the main point is that we have to avoid in the proofs the concept of entropy and the variational principle of pressure.

\bigskip
Remember that we denote by $ \mathcal M_\sigma$ the set of $\sigma$ invariant Borel probability measures over ${\cal  B}$. As $ \mathcal M_\sigma$ is compact, given $A$, there always exists a subsequence $\beta_n$, such that $\mu_{\beta_n A}$ converges to an invariant probability measure.

We  consider the following problem: given $A:{\cal B} \to \mathbb{R}$ Lipschitz, we want to find measures that  maximize, over $ \mathcal M_\sigma$, the value

$$ \int A(x) \,d \mu(\mathbf {x}). $$
\bigskip

We define  
$$m(A)=\max_{\mu\in\mathcal M_\sigma} \left\{ \int A d\mu \right\}\,.$$

Any of these measures will be called a maximizing probability measure, which is sometimes denoted by $\mu_\infty$. As $\mathcal M_\sigma$ is compact, there exist always at least one maximizing probability measure. It is also true that there exists ergodic maximizing probability measures. Indeed, the set of maximizing probability measures is convex, compact and the extreme probability measures of this convex set are ergodic (can not be expressed as convex combination of others \cite{Kel}).  Any maximizing probability measure is a convex combination of ergodic ones \cite{PY}.

Even when $A$ is Holder the maximizing probability measure $\mu_\infty$  do not have to be unique. For instance, suppose that $A$ is Holder and has maximum value just in the union of two different fixed points (for the shift $\sigma$) $p_0\in {\cal B}$ and $p_1\in {\cal B}$ . In this case the set of maximizing probability measures $\mu_\infty$ is $\{ t \, \delta_{p_0} + (1-t) \delta_{p_1}|\,t\in [0,1]\}$.

Note that $\delta_{p_0}$ and $\delta_{p_1}$ are ergodic, but the other maximizing probability measures are not.

Similar definitions for a potential $A:{\cal B}_i \to \mathbb{R}$ and maximization of $\int A \,d \hat{\mu}$,  over all the $\hat{\mu}$ which are $\hat{\sigma}$-invariant probability measures, can also be considered. Questions about selection of measure also make sense.

\begin{definition}\label{sub}
A continuous function $u: {\cal B} \to \mathbb{R} $ is called a
{\em calibrated subaction} for $A:{\cal B}\to \mathbb{R}$, if, for any $y\in {\cal B}$, we have

\begin{equation}\label{c} u(y)=\max_{\sigma(x)=y} [A(x)+ u(x)-m(A)].\end{equation}

\end{definition}

This can also be expressed as
$$m(A)= \max_{a \in M} \{A(ay)+  u (ay) - u(y) \} .$$

Note that for any $x\in {\cal B}$ we have
 $$ u(\sigma(x)) -  u(x) - A(x) + m(A) \geq0.$$

The above equation for $u$ can be seen as a kind of discrete version of a sub-solution of the Hamilton-Jacobi equation \cite{CI} \cite{BC} \cite{Fathi}. It can be also seen as a kind of dynamic additive eigenvalue problem \cite{CD} \cite{CG} \cite{GL3}.

If $u$ is a calibrated subaction, then $u+c$, where $c$ is a constant, is also a calibrated subaction. An interesting question is when such calibrated subaction $u$ is unique up to an additive constant.

Remember that if $\nu$ is invariant for $\sigma$, then for any continuous function $u:{\cal B}\to \mathbb{R}$ we have
$$ \int \,[ u(\sigma(x)) -  u(x)]\, d \nu=0$$

Suppose $\mu$ is maximizing for $A$ and $u$ is a calibrated subaction for $A$.

 It follows at once (see for instance \cite{CLT} \cite{J1} \cite{Sou} for a similar result) that for any $x$ in the support of $\mu_\infty$ we have
\begin{equation}\label{funcaoR} u(\sigma(x)) -  u(x) - A(x) + m(A)=0. \end{equation}

 In this way if we know the value $m(A)$, then a calibrated subaction $u$ for $A$ helps us to identify the support of maximizing probabilities. The above equation can be true outside the union of the supports of the  maximizing probabilities.

 Maximizing probability measures 
  are natural candidates for being selected by $\mu_{\beta A}, $ as $\beta \to \infty$.
 But, in our setting, without the maximizing principle of pressure (which one can take advantage of the classical Thermodynamic Formalism) this is not so obvious. We address the question in section \ref{sec2b}.

\begin{proposition} For any $\beta$, we have
$-\|A\| < \frac{1}{\beta}\log\lambda_{\beta} < \|A\|$.

\end{proposition}

 {\it Proof: }
 Fix $\beta>0$. We choose $\bar x$ the maximum of $\psi_{\beta A}$ in $\cal{ B}$ and $\tilde x$ the minimum of $\psi_{\beta A}$ in $\cal{ B}$.
  Now, if $\|A\|$ is the uniform norm of $A$, we have
 $$\lambda_{\beta}=\frac{1}{\psi_{\beta A}(\bar x)} \int e^{\beta A(a\,\bar x)} \psi_{\beta A}(a\, \bar x)d a\leq
 \int e^{\beta A(a\,\bar x)}d a
 \leq e^{\beta \|A\|} \;\mbox{and}$$
 $$\lambda_{\beta}=\frac{1}{\psi_{\beta A}(\tilde x)} \int e^{\beta A(a\,\tilde x)} \psi_{\beta A}(a \,\tilde x)da \geq
  \int e^{\beta A(a\,\bar x)}d a
 \geq e^{-\beta \|A\|} \; ,$$
which proves the result.
\cqd

\vspace{0.3cm}
From now on, we will suppose $M=\mathbb{S}^1$ to avoid technical issues. But we claim that the following results hold for more general connected and compact manifolds.

Considering a subsequence $\beta_n$ we get the existence of a limit $ \frac{1}{\beta_n}\log\lambda_{\beta_n}\to K$, when $n \to \infty$.
 By taking a subsequence we can assume that is also true that there exists $V$ Lipschitz, such that
$V:=\lim_{n\to\infty}
\frac{1}{\beta_n}\log(\psi_{\beta_n A}).$

Given $y \in \cal{B}$, consider the equation
 $$\lambda_{\beta_n}=\frac{1}{\psi_{\beta_n A}(y)} \int e^{\beta_n A(a\, y)} \psi_{\beta_n A}(a\,y)d a .$$

 It follows from Laplace method that, when $\beta\to \infty$,
 $$K= \max_{a \in \mathbb{S}^1} \{A(ay)+  V (ay) - V(y) \} .$$

If we are able to show that $K=m(A)$, then we can say that any limit of subsequence $\lim_{n\to\infty}
\frac{1}{\beta_n}\log(\psi_{\beta_n A})$ is a calibrated subaction, and we will get, finally, that
$$\lim_{\beta \to \infty}  \frac{1}{\beta}\log\lambda_{\beta\, A}=\lim_{\beta \to \infty}  \frac{1}{\beta}\log\lambda_{\beta}= m(A).$$

Next theorem is inspired by Theorem 1 in \cite{Bousch-walters} and Theorem 3.3 in \cite{Gomes1}.
It follows from the last part of its proof  that $K=m(A)$.

\begin{theorem}\label{existsubac} Given $A$ Lipschitz there exists $u$ Lipschitz which is a calibrated subaction for $A$. As a consequence, we have that
$$\lim_{\beta \to \infty}  \frac{1}{\beta}\log\lambda_{\beta}= m(A).$$
\end{theorem}

\noindent{\bf Proof.} Suppose $A:{\cal B} \to \mathbb{R}$ is Lipschitz.

Given $0<\lambda\leq 1$, consider the operator $\hat{{\cal L}}_\lambda:{\cal C} \to {\cal C}$ given by,

$$ \hat{{\cal L}}_\lambda(u)(x) = \sup_{a\in \mathbb{S}^1} [A(ax) + \lambda u(ax)].$$

Given  $x\in {\cal B}$ , we denote by $a_x\in \mathbb{S}^1$ one of the points $a$ where the supremum is attained.

It is easy to see that for any $0<\lambda<1$, the transformation $\hat{{\cal L}}_\lambda$
is a contraction on ${\cal C}$ with the uniform norm.
Indeed, given $x \in {\cal B}$
$$  \sup_{a\in \mathbb{S}^1} [A(ax) + \lambda u(ax)]-  \sup_{b\in \mathbb{S}^1} [A(bx) + \lambda v(bx)]\leq $$
$$ [A(a_x\,x) + \lambda u(a_x\,x)]- [A(a_x\,x) + \lambda v(a_x\,x)]\leq$$
$$ \lambda u(a_x\,x)- \lambda v(a_x\,x)\leq\lambda ||u-v||.$$

Denote by $u_\lambda$ the corresponding fixed point in ${\cal C}$.
We want to show that $u_\lambda $ is equicontinuous.
Consider $x_0,y_0\in{\cal B}$. For the given $x_0$ we take the corresponding $a_{x_0}\in M$, and then the we get $x_1=a_{x_0}x_0.$ By induction, given $x_j$, get $ x_{j+1}=a_{x_j}x_j.$

We can also can get a sequence $y_j\in {\cal B}$, $j \geq 1$, such that, $ y_j = a_{x_{j-1}\,...\, } a_{x_1} a_{x_0} y_0$. Note that for all $j$ we have $\sigma^j (y_j)=y_0.$

As for any $j$ we have $u_\lambda(y_j) \geq A(y_{j+1}) - \lambda u_\lambda (y_{j+1})$, then
$$ u_\lambda (x_j) - u_\lambda (y_{j})\leq $$
$$  [A(x_{j+1})- A(y_{j+1})] + \lambda \, [ u_\lambda (x_{j+1}) - u_\lambda (y_{j+1})].$$

Therefore, given $x_0,y_0$
$$ u_\lambda (x_0) - u_\lambda (y_{0})\leq \sum_{j=0}^\infty \lambda^j \,[ A (x_{j}) - A (y_{j})]\leq$$
$$ (1-\lambda) \sum_{j=0}^\infty \lambda^j \sum_{i=0}^j\,[ A(x_i) - A(y_i)]\leq$$
$$ \sup_j \sum_{i=0}^j\,[ A(x_i) - A(y_i)]\leq$$
$$ ||A|| \,  \sup_j \sum_{i=0}^j\,(\frac{1}{2})^j d(x_0,y_0)< ||A|| \,2 \,d(x_0,y_0).$$

This shows that $u_\lambda$ is Lipschitz, and, moreover, that $u_\lambda$, $0\leq \lambda<1$, is an equicontinuous family. Note the very important point: the Lipschitz constant of $u_\lambda$ depends on $||A||$.


Denote $u^*_\lambda = u_\lambda - \max u_\lambda$. Using Arzela-Ascoli we get the existence of a subsequence $\lambda_n\to 1$ such that
$u^*_{\lambda_n }\to u$.

We claim that $u$ is a subaction.

Indeed, given $x\in {\cal B}$,
as $|u_\lambda(x)| \leq \lambda \, |u_\lambda(a_x\,x) | + |A(a_x\,x)|\leq \lambda\, ||u_\lambda|| + ||A(x)||,$
then $(1-\lambda) ||u_\lambda||< C$, where $C$ is a constant.

From this follows that there  is a constant $k$, such for some subsequence (of the previous subsequence $\lambda_n$), which will be also denoted by $\lambda_n$, we have
$(1-\lambda_n) ||u_{\lambda_n}||\to k$.

Note that for any $\lambda$
$$u^*_\lambda (x)= u_\lambda(x) - \max u_\lambda= $$
$$-(1-\lambda) \max u_\lambda + u_\lambda (x) - \lambda \max u_\lambda=$$
$$-(1-\lambda) \max u_\lambda + \max_{a \in \mathbb{S}^1} \{A(ax)+ (\lambda  u_\lambda (ax) - \lambda \max u_\lambda)\}.$$

Taking the limit $n$ to infinity for the sequence $\lambda_n$ we get
$$u(x) =-k +\max_{a \in \mathbb{S}^1} \{A(ax)+  u (ax) \}= \max_{a \in \mathbb{S}^1} \{A(ax)+  u (ax) -k \} .$$

Now, all we have to show is that $k=m(A).$

From the above it follows at once that

$$ -u(\sigma(y)) + u(y) + A(y)\leq k.$$

If $\nu$ is a $\sigma$-invariant  probability measure, then,
$$\int A(y) d \nu(y)=\int  [u(\sigma(y)) - u(y) + A(y)]\,d \nu(y)\leq k,$$
and, this shows that $m(A)\leq k.$

Now we show that $m(A) \geq k$. Note that for any $x$ there exist $y=a_x\,x$ such that $\sigma(y)=x$, and
$$ -u(\sigma(y)) + u(y) + A(y)= k.$$

Therefore, the compact set $K= \{y \,|\,-u(\sigma(y)) + u(y) + A(y)= k\}$ is such that, $K' = \cap_n \, \sigma^{-n} (K)$ is non-empty, compact and $\sigma$-invariant. If we consider an $\sigma$-invariant probability measure $\nu$ with support on $K'$, we have that $\int A(y) d \nu(y)=k$. From this follows that $m(A) \geq k$.

\cqd

\vspace{0.3cm}


Now we state a general result assuming just that $A$ is continuous (not necessarily Lipschitz). We refer the reader to Theorem 1 in \cite{GL1}, Proposition 4 in \cite{LMST}, Theorem 2.4 in \cite{Gomes1} for  related results.

\begin{theorem}\label{formuladual}
Given a potential $ A \in {\cal C} $, we have
$$ m(A) =
\inf_{f \in {\cal C}} \,\,\,\max_{(\mathbf a, \mathbf x) \,\in\,  \mathbb{S}^1 \times {\cal B}} \,\,\,[A(\mathbf a\, \mathbf x) + f(a\mathbf x)  - f(\mathbf x))]. $$
\end{theorem}

{\it Proof:} First, consider the convex correspondence
$ F: {\cal C} \to \mathbb R $ defined by $ F(g) = \max (A + g) $.
Consider also the subset
$$ \mathcal G = \{g \in {\cal C}: \text{there exists} \,\,f \,\,\text{such that} \,\,g(\mathbf a \mathbf x) =
f(a\, \mathbf x) - f(\mathbf x),\,\, \,\,f \in {\cal C} \}\neq \emptyset. $$

 Now consider the concave correspondence $ G: {\cal C} \to \mathbb R \cup \{- \infty \} $ taking
$ G(g) = 0 $, if $ g \in \bar{\mathcal G} $, and $ G(g) = - \infty $ otherwise.

Let $ \mathcal S $ be the set of the signed measures over the Borel
sigma-algebra of $ {\cal B} $. Remember that the corresponding Fenchel transforms,
$ F^*: \mathcal S \to \mathbb R \cup \{+ \infty \} $ and $ G^*: \mathcal S \to \mathbb R \cup \{- \infty \} $, are given by
$$ F^* (\hat \mu) = \sup_{g \in {\cal C}} \left [\int_{} g(\mathbf a \mathbf x) \; d\hat\mu(\mathbf a \mathbf x) - F(g) \right] \;
\text {, \,and} $$
$$ G^* (\hat \mu) = \inf_{g \in {\cal C}} \left [\int g(\mathbf a \mathbf x) \; d\hat\mu(\mathbf a \mathbf x) - G(g) \right]. $$
Denote
$$ \mathcal S_0 = \left \{\hat \mu \in \mathcal S:  \int f(\mathbf a\mathbf x) \; d\hat\mu(\mathbf a \mathbf x) =
\int f(\mathbf x) \; d\hat\mu(\mathbf x) \; ,\; \forall \; f \in {\cal C} \right \}. $$

We denote by ${\cal M}$ the set of probability measures over ${\cal B}$.

Given $ F $ and $ G $ as above, we claim that
$$ F^* (\hat \mu) =
\left\{ \begin{array}{ll} {\displaystyle - \int_{\hat \Sigma} A(\mathbf y, \mathbf x) \; d\hat\mu(\mathbf y, \mathbf x)}
& \mbox{if $ \hat \mu \in \mathcal M $} \\ + \infty & \mbox{otherwise} \end{array} \right.
\; \text{ and} $$
$$ G^* (\hat \mu) =
\left\{ \begin{array}{ll} 0 & \mbox{if $ \hat \mu \in \mathcal S_0 $} \\ - \infty & \mbox{otherwise} \end{array} \right.. $$

We refer the reader to the \cite{GL1} or \cite{LMST} for a proof of this claim (which is basically the same as we need here).

Once the correspondence $ F $ is Lipschitz, the theorem
of duality of Fenchel-Rockafellar \cite{Rock} assures
$$ \sup_{g \in {\cal C}} \left[G(g) - F(g)\right] = \inf_{\hat \mu \in \mathcal S} \left[F^*(\hat \mu) - G^*(\hat \mu)\right]. $$


$$ \sup_{g \in \mathcal G} \left [- \max_{(\mathbf a, \mathbf x) \in \mathbb{S}^1 \times {\cal B}} (A + g) (\mathbf a \mathbf x) \right] =
\inf _{\hat \mu \in \mathcal M_\sigma} \left [- \int A(\mathbf a \mathbf x) \; d\hat\mu(\mathbf a \mathbf x) \right]. $$
Finally, from the definition of $ \mathcal G $, the claim of the theorem follows.

\cqd

\vspace{0.6cm}

\section{A definition of entropy for Gibbs states at positive temperature and selection of probability measure}\label{sec2b}

Given a  Lipschitz function $A$ we have that

 $$\int\, \frac{e^{A(ax)}\evc(ax)}{\evl \evc(x)}da=1\,,\,\forall x \in \mathcal{B}\,.$$

 We denote as before $$\bar A = A+\log  \evc-\log \evc \circ \sigma -\log \evl,$$
where $\sigma:\BB \rar \BB$ is the usual shift map.
In this case the normalized potential $\bar A$ satisfies
$$\II e^{\bar A (ax)}da = 1 \,,\,\forall x \in \mathcal{B}\,,$$
 which means $\mathcal{L}_{\bar A}(1)=1$.

 Therefore,
 $$ \int_{\mathcal{B}}\,\left[ \int_M\, e^{\bar A (ax)}\, da\right] \, d \, \mu_A(x) =1.$$

Note that for a fixed $x$ the value $\bar A(ax)$ can not be smaller than zero for all $a\in M$.
This is quite different from the analogous case where we consider the shift over $\{1,2..,d\}^\mathbb{N}$ in the classical Thermodynamic Formalism.

For each $a\in M, \,x\in {\cal B}$, we denote by $J(ax) = \min \{1,  e^{\bar A (ax)}\}.$

\begin{definition}
Given the invariant probability measure $\mu_A$, associated to the Lipschitz potential $A$, we define the entropy of $\mu_A$ as
$$ h(\mu_A) = -\int\,   \, \log J(y)\,  d\mu_A (y)>0 .$$

\end{definition}

In other words
$$ h(\mu_{A}) =-\, \int \bar A(y)\, \,I_{\{ \bar A \leq 0\}}\,(y)\,\,\, d \mu_{ A}(y).$$

The set of probability measures $\mu_A$, with $A$ Lipschitz, is dense in the set of $\sigma$-invariant probability measures \cite{L3}.

Note that $\mu_A$ is $\sigma$-invariant
$$ - h(\mu_A)= \int\,  \, \log J(y)\,   d\mu_A (y) \leq$$
$$\int\,\, \log\left(\,\frac{e^{A(y)}\evc(y)}{\evl \evc( \sigma(y) )}\,\right)\,d \mu_A (y)= \int A d \mu_A - \log \lambda_A.$$

Therefore,
$$ \log \lambda_A \leq h(\mu_A) + \int A d \mu_A.$$

For a fixed $A$ consider now for each real value $\beta$ the corresponding potential $\beta A$. Therefore,
$$ \log \lambda_{\beta A} \leq h(\mu_{\beta A}) + \beta \,\int A\, d \,\mu_{\beta A}.$$

Suppose for a certain subsequence $\beta_n$ we have that $\mu_{\beta_n \,A} \to \mu$.

If we divide the last inequality by  $\beta_n$, and, taking limit in $n$, we get
$$ m(A) \leq \limsup_{n \to \infty} \frac{h(\mu_{\beta_n A})}{\beta_n}  + \int A \,d \mu.$$

From the above we can derive:

\begin{theorem}

Suppose that  $\mu= \lim_{n\to \infty} \mu_{\beta_n A}$, for some subsequence $\beta_n$, and
$$\limsup_{n \to \infty} \frac{h(\mu_{\beta_n A})}{\beta_n}=0,$$
then, the limit measure $\mu$ is a maximizing probability measure.

\end{theorem}

\begin{corollary}
 If the maximizing probability measure $\mu_\infty$ for $A$ is unique,
 and,
 $$\limsup_{\beta \to \infty} \frac{h(\mu_{\beta A})}{\beta}=0,$$
then, $\mu_{\beta A}$, when $\beta \to \infty$, selects the maximizing probability measure $\mu_\infty$.

\end{corollary}

\vspace{0.6cm}

\section{Analysis of the case in which the potential depends on two coordinates}\label{sec3}

In this section we suppose the potential depends on two coordinates and the metric space is $M=\mathbb{S}^1$.
In this case the Ruelle operator has a simple form. We will make the usual identification of $\mathbb{S}^1$ with $[0,1]$  (in further sections we will make the identification of $\mathbb{S}^1$ with $[0,2\pi]$). We will present several results from \cite{LMST} which will be needed in future sections.

We will need to define the following operators:

\begin{definition} Let
${\cal L}_{\beta},{\bar{\cal L}}_{\beta}:C([0,1])\to C([0,1])$ be given by

\begin{equation}\label{L}   {\cal L}_{\beta} \psi (y) =\int e^{\beta A(x,y)} \, \psi(x)dx,  \end{equation}

\begin{equation}\label{bar L}   {\bar{\cal L}}_{\beta} \psi (x) =\int e^{\beta A(x,y)} \, \psi(y)dy.  \end{equation}
\end{definition}
\vspace{0.2cm}
We refer the reader to \cite{Ka} and \cite{Sch} chapter IV for general results on positive integral operators. The next theorem (Krein-Ruthman) is well known.
It will follow  that,
when $A$ depends just on two coordinates $(x_0,x_1)$, then the eigenfunction of the Ruelle operator (as defined in previous sections) depends only on the first coordinate $x_0$ (similar to \cite{Sp}).

\begin{theorem}\label{KreinRuthman} 
The operators  ${\cal L}_{\beta}$ and ${\bar{\cal L}}_{\beta}$ have the same  positive maximal eigenvalue $\lambda_{\beta}$, which is simple and isolated. The eigenfunctions associated are  positive functions.
\end{theorem}

Let us call $\psi_{\beta},\bar{\psi}_{\beta}$  the
positive eigenfunctions for ${\cal L}_{\beta}$ and ${\bar{\cal
L}}_{\beta}$ associated to $\lambda_{\beta}$, which satisfy the normalization condition $\int
\psi_{\beta}(x)\,dx=1$ and $\int \bar\psi_{\beta}(x)\,dx=1$.

\vspace{0.3cm}

We will define a density  $\theta_{\beta}:[0,1]\to \mathbb{R}$ by \begin{equation}\label{theta}\theta_{\beta}(x):=
\frac{\psi_{\beta}(x)\,\, \bar\psi_{\beta}(x)}{\pi_{\beta}},\end{equation} where $\pi_{\beta}=\int \psi_{\beta}(x) \bar\psi_{\beta}(x)dx$,
and a transition $ K_{\beta}:[0,1]^2\to \mathbb{R}$ by
\begin{equation}\label{K}K_{\beta}(x,y):=\frac{e^{\beta A(x,y)}\,\,\bar \psi_{\beta}(y)}{\bar \psi_{\beta}(x)\,\lambda_{\beta}}\,. \;\end{equation}

The above expressions are consistent with the results obtained in \cite{FH} section 3. This can be formulated also as a variational pressure problem as we will see soon.

Note that if $A(x,y)=A(y,x)$, then $\psi_{\beta}$ and $\bar \psi_{\beta}$ are constant, and, therefore $\theta_\beta$ is constant equal to $1$. This happen for $M=\mathbb{S}$ in the case $A$ is of the form $U(x-y)$ for a periodic function $U$. This case will be consider later.

Consider a probability measure $\nu$ on $ [0,1]^2$ that can be disintegrated as $d\nu(x,y)= d\theta(x) dK_x(y)$, where $\theta:[0,1]\to [0,+\infty)$ and $K:[0,1]^2\to [0,+\infty)$ are continuous functions.  We will denote this by $\nu=\theta K$, where
$\theta$ is a continuous density of probability on $[0,1]$.
\begin{definition} \label{est} A   probability measure $\theta$ on $[0,1]$ is called
{\em stationary} for a transition  $K(\cdot,\cdot)$, if
$$\theta (B)=\int K(x,B)d\theta(x),\;\;\;\;\;\;\;\mbox{ for all interval}\,\, B\in [0,1]\,.$$

\end{definition}

More explicitly we assume
$K:[0,1]^2\to[0,+\infty)$ and  $\theta:[0,1]\to[0,+\infty)$ satisfy the following equations:

\begin{equation}\label{transicao}
\int K(x,y)\, dy=1, \hspace{1cm} \forall \, x \in [0,1],
\end{equation}

\begin{equation}\label{normalizada}   \int \theta(x)\,K(x,y) \,dx dy=1, \end{equation}

\begin{equation}\label{equilibrio}  \int \theta(x)\, K(x,y) \,dx = \theta(y),\hspace{1cm} \forall \, y \in [0,1].
\end{equation}

Given the initial probability measure $\theta$  and the transition $K$, as
above, one can define a Markov process $\{X_n\}_{n\in \mathbb{N}}$
with state space $[0,1]$ (see \cite{LMST} for more details). The measure $\mu$ over $[0,1]^\mathbb{N}$ which describes this process is

$$ \mu(A_0...A_n\times [0,1]^\mathbb{N}):=\int_{A_0...A_n}\,\theta(x_0)\,
K(x_{0},x_1)...K(x_{n-1},x_n)\,dx_n...dx_0 $$
for any cylinder $ A_0...A_n
\times [0,1]^\mathbb{N}$.


If $\theta$ is stationary the Markov Process $X_n$ will be stationary. 

Note that $\theta_\beta$ above is stationary for $K_\beta(x,y)$. In this way we  can define  $\nu_\beta= \theta_\beta K_\beta  \,$ on $[0,1]^2$.


For instance,
$$\mu_{\beta,A}( \,[a_1,a_2]\times [b_1,b_2]\times [c_1,c_2]\,\times [0,1]^\mathbb{N})=$$
\begin{equation}\label{medidaspitzer}
     =\int_{a_1}^{a_2} \int_{b_1}^{b_2} \int_{c_1}^{c_2}
\theta_\beta (x_0)\, K_\beta(x_0,x_1)\,K_\beta(x_1,x_2)dx_2\, d x_1 \,dx_0.
\end{equation}

The next result is similar to the one described in \cite{Sp}.

\begin{theorem}\label{spitzer_vs_PP}
 Suppose $A$ is a Holder continuous function.  Then the probability measure $\mu_{\beta, A}$ defined in \eqref{medidaspitzer} is the Gibbs state for the potential $\beta A$.
\end{theorem}

\noindent{\bf Proof.}
 We need to show that $\mathcal{L}_{\bar A}^*(\mu_{\beta, A})=\mu_{\beta, A}$, where $\beta\bar A = \beta A+\log  \psi_{\beta }-\log \psi_{\beta }\circ \sigma -\log \lambda_{\beta }$. Indeed, let $g \in \mathcal{C}$ such that $g(x_0,x_1,...)=g(x_0,...,x_k)$, by definition of $\mathcal{L}_{\bar A}^*$ we have
$$  \int_{\BB} g\, d\mathcal{L}^*_{\bar A}(\mu_{\beta, A}) =  \int_{\BB} \mathcal{L}_{\bar A}(g)\, d\mu_{\beta, A}\,$$
$$=\int_{[0,1]^{k}}\,\bigg[ \int_{[0,1]}\,e^{\beta\bar A(a,x_0)}g(a,x_0,...,x_{k-1})\, d\,a\bigg]\,\theta_\beta (x_0)\prod_{j=0}^{k-2} K_\beta(x_{j},x_{j+1}) dx_{k-1}...dx_0$$
$$=\int_{[0,1]^{k}}\bigg[ \int_{[0,1]}\,e^{\beta A(a,x_0)}\frac{ \psi_\beta (a)}{\lambda_{\beta }\psi_\beta(x_0)}\, g(a,...,x_{k-1})\, d\,a\bigg]\,\theta_\beta (x_0)\prod_{j=0}^{k-2} K_\beta(x_{j},x_{j+1})dx_{k-1}...dx_0$$

$$=\int_{[0,1]^{k+1}}\,\,e^{\beta A(a,x_0)}\, \frac{ \psi_\beta (a)}{\lambda_{\beta}}\, g(a,...,x_{k-1})\, \,\frac{\bar\psi_\beta (x_0)}{\pi_\beta}\prod_{j=0}^{k-2} K_\beta(x_{j},x_{j+1})\,dx_{k-1}...dx_0\,da$$

$$=\int_{[0,1]^{k+1}}\,\, g(a,...,x_{k-1})\, e^{\beta A(a,x_0)}\frac{\bar\psi_\beta(x_0)}{\lambda_{\beta }\bar\psi_\beta(a)}\,\,\frac{\bar\psi_\beta (a)\,\psi_\beta(a)}{\pi_\beta } \,\prod_{j=0}^{k-2} K_\beta(x_{j},x_{j+1})\,\,dx_{k-1}...dx_o da$$

$$=\int_{[0,1]^{k+1}}\,\, g(a,x_0,...,x_{k-1})\,\theta_{\beta}(a) \,K_\beta(a,x_{0}) \prod_{j=0}^{k-2} K_\beta(x_{j},x_{j+1})\, \, \,dx_{k-1}...dx_0 da$$

$$=\int_{[0,1]^{k+1}}\,\, g(x_0,x_1,...,x_{k-1}, x_{k})\,\theta_{\beta}(x_0)  \,\prod_{j=0}^{k-2} K_\beta(x_{j},x_{j+1}) K_\beta(x_{k-1},x_{k}) \,dx_{k}...dx_1\,dx_0.$$

Hence, for any continuous $g$
$$\int_{\BB} g(x_0,...,x_k)\, d\mathcal{L}^*_{\bar A}(\mu_{\beta, A}) =\int_{\BB} g(x_0,...,x_k) d\mu_{\beta, A}.$$
\cqd


The entropy (as defined in section \ref{sec2b}) of such probability measure $\mu_{\beta A}$ is
$$ h(\mu_{\beta A}) =-\, \int A(y)\, \,I_{\{ A- \log \lambda_\beta\leq 0\}}\,(y)\,\,\, d \mu_{\beta A}(y) + \log \lambda_{\beta}.$$

\begin{definition}
We denote by ${\cal M}_0$ the set of all $\nu=\theta K$ on $[0,1]^2$, where $\theta$  is stationary for $K$.
\end{definition}

\vspace{.3cm}
\begin{definition}
For an absolutely continuous
probability measure $\nu\in\mathcal M_{[0,1]^2}$, given by a density $\nu(x,y)dxdy$, we denote $S[\nu]$ by

\begin{equation}
S[\nu]=-\int \nu(x,y) \log\left(   \frac{\nu(x,y)}{\int\nu(x,z)dz}      \right) dx dy \,.
\end{equation}
\end{definition}
\noindent{\bf Remark.}

The $S[\nu]$ was called "penalized entropy" in \cite{LMST} and it is a kind of relative entropy with respect to Lebesgue measure. It is different kind of definition of entropy  (from the previous one we consider before).

It is easy to see that  any $\nu=\theta K\in {\cal M}_{0}$  satisfies
\begin{equation}
S[\theta P]=-\int \theta(x) K(x,y) \log\left(    K(x,y)      \right) dx dy \,.
\end{equation}

The value $S[\theta\, K]$ assume negative values.

We can consider now the variational problem

\begin{equation}\label{formterm}   P(A)=\max_{\nu=\theta K\in{\cal M}_{0}} \left\{ \int\beta  A(x,y)\, d\nu +
 S[\nu]\right\}.
\end{equation}

This is equivalent to maximize
$$\max_{\nu=\theta K\in{\cal M}_{0}} \left\{ \int\beta  A(x,y)\theta(x) K(x,y) dx dy -
 \int \theta(x) K(x,y) \log\left(    K(x,y)      \right) dx dy \right\}
$$

\begin{definition}

A probability measure $\nu$ in ${\cal M}_0$ is called an equilibrium state for $A$ (which depends on two coordinates) if attains the maximal value $P(A)$. The value $P(A)$ is called the pressure (or Free Energy) of $A$
\end{definition}


We refer the reader to \cite{LMST} for the proof of the following result.

\begin{proposition}\label{maxMarkov} The stationary measure $\nu_{\beta}=\theta_{\beta}K_{\beta}$ defined above maximize
 $$\beta\,  \int    A(x,y)\;d\nu+ S[\nu],$$
 over all stationary $\nu=\theta K\in {\cal M}_0$. Also

 $$P(A)=\log \lambda_{\beta}=\int \beta A \;\theta_{\beta} K_{\beta}dxdy +S[\theta_{\beta} K_{\beta}].$$
\end{proposition}

When the potential $A$ depends just on two coordinates the equation used in the definition of subaction can be simplified.

\begin{definition}\label{sub}
A continuous function $u: [0,1] \to \mathbb{R} $ is called a $[0,1]$-
{\em calibrated forward-subaction}  if, for any $y\in [0,1]$, we have

\begin{equation}\label{c} u(y)=\max_{a\in [0,1]} [A(ay)+ u(a)-m(A)].\end{equation}

\end{definition}

We refer the reader to \cite{CG} for related problems in a different setting. The equation for $u$ above also appears in problems related to the additive eigenvalue \cite{CD} \cite{CG}.

A function $u$ as above can be seen as a function on $x \in [0,1]^\mathbb{N},$ where $x=(x_0,x_1,x_2,x_3,...)$, which depends just on the first coordinate $x_0$. Therefore, a $[0,1]$-calibrated forward-subaction is a also
calibrated subaction (in the previous sense). We point out that $[0,1]$-
calibrated forward-subactions do exist (see \cite{LMST}).

 An interesting question on the case of selection of measures $\mu_\beta\to \mu_\infty$ is: what happens with the measure of a
 particular subset $D$ of ${\cal B}$ when $T\to 0$ (or, $\beta \to \infty$)? A Large Deviation Principle (see \cite{DZ} for general references)  is true under certain conditions. We refer the reader to \cite{LMST} for the proof of the result below.

\begin{theorem}\label{teoLDP}

 If $A$  has only one   maximizing probability measure $\mu_\infty$ and there exist an unique  $[0,1]$- calibrated  forward-subaction $V$ for $A$, then the following LDP is true:
for each cylinder $D=A_0....A_k\times [0,1]^\mathbb{N}$, the following limit exists

$$\lim_{\beta\to\infty}\frac{1}{\beta}\ln
\mu_{\beta A}(D)= -\inf_{\mathbf{x}\in D}I(\mathbf{x})\,.$$

where
$I:[0,1]^{\mathbb{N}}\to [0,+\infty]$ is a function defined by
$$I(\mathbf{x}):= \sum_{i\geq 0}V(x_{i+1})-V(x_i)-(A-m(A))(x_i,x_{i+1})\,.   $$

\end{theorem}

Results about Large deviations in the setting of Thermodynamic Formalism appear in \cite{CLT} \cite{LM}.

\begin{definition} We say that $A:[0,1]^2 \to \mathbb{R}$ satisfies the twist condition, if $A$ is $C^2$, and
$$\frac{\partial^2 A}{\partial x \partial y} \neq 0.$$
\end{definition}

This property is an open condition under the right topology.

The next theorem (see \cite{LMST} for a proof) addresses the question of uniqueness  when we add a magnetic term $f(x)$ to $A(x,y).$ Related results in a different setting appear in \cite{Ban} \cite{BC}. The above condition for $A$ replaces the convexity of the Lagrangian which is crucial in Aubry-Mather theory \cite{CI}.

\begin{definition} We will say that a property  is {\em generic    for  $A$, $A\in C^{2}([0,1]^2)$, in Ma\~n\'e's sense}, if
the property is true for $A+f$, for any $f$, $f\in C^{2}([0,1])$, in a set $G$ which is generic (in Baire sense).
\end{definition}

This concept was initially  introduced in the Aubry-Mather  setting in \cite{Mane}.

We will show below that under the twist condition the uniqueness of $[0,1]$-forward backward-subaction is generic in Ma\~n\'e's sense.

\begin{theorem} \label{principal1} Consider the class of all $A:[0,1]^2 \to \mathbb{R}$
which is  $C^2$ and satisfies the twist condition
$\frac{\partial^2 A}{\partial x \partial y} \neq 0$, then
 there exists a generic set $\mathcal{O}$  in $C^2([0,1])$  (in Baire sense) such that:
 \medskip

 (a) for each $f \in \mathcal O$, $f:[0,1]\to \mathbb{R}$, for "any" $A$ we have that
given $\mu, \tilde \mu\in {\cal M}_\sigma$ two maximizing measures for $A+f$ (i.e., $m(A+f)=\int (A+f)\;d\mu=\int (A+f)\;d\tilde\mu$),
then
$$\nu=\tilde \nu,$$
where $\nu$ and $\tilde \nu$ are the projections of $\mu$ and $\tilde \mu$ in the first two coordinates.

\bigskip

(b) for "any" $A$ the  $[0,1]$-calibrated forward-subaction  for  $A+f$ is unique, for each $f \in \mathcal O$ (up to an additive constant).
\end{theorem}

In the above theorem the potential $A$ is considered the interaction and $f $ the magnetic term.
Therefore, it claims, among other things,  that for "any" $A$ we have uniqueness of the calibrated subaction (up to an additive constant) for a generic magnetic term $f$.

The next theorem (see \cite{LMST} for a proof) addresses the question of the graph property for a probability measure. A related result in the setting of Thermodynamic formalism appears in \cite{LOT}.

\begin{theorem}
\label{principal2} If $A:[0,1]^2 \to \mathbb{R}$
is  $C^2$, and satisfies the twist condition
$\frac{\partial^2 A}{\partial x \partial y} \neq 0$, then,
the projected measure $\nu$ on $[0,1]^2$ of the maximizing probability measure $\mu_\infty$ (on ${\cal B}$) has support on a graph.
\end{theorem}

The problem we consider above can be seen as a Transshipment Problem (see \cite{LMST}).
For related results see also \cite{GL3} and \cite{CLO}.

The graph property of a measure is of great importance in Aubry-Mather Theory \cite{CI} \cite{Fathi} \cite{Mat}.
\bigskip

\vspace{0.6cm}

\section{DLR Gibbs Measures and Transfer Operator}\label{sec-dic}
\bigskip

\subsection{One-Dimensional Systems and Transfer Operator}

Given the potential $A$
we will use the following terminology: Gibbs-TF for $A$ denotes the set of measures usually considered in the Thermodynamical Formalism
(as, for example, in \cite{PP}, or in the first part of this paper) and Gibbs-DLR for $A$
the set of measures constructed as in the Dobrushin-Lanford-Ruelle formulation of Statistical Mechanics,
where the Gibbs measures are obtained from Specification Theory point of view, for a complete exposition see \cite{G,Ruelle,Si,Sa}.
For reasons that will be clarified latter we adopt the notation $\mu^{ A,\sigma'}$ for this measures,
where $\sigma'$ is an element of the state space which is called sometimes a boundary condition.

The measures obtained by the the first construction (Section \ref{sec1}) are denoted here by $m=m_A$, and they are
defined over the $\sigma$-algebra of
$\mathcal{B}=(S^1)^{\mathbb{N}}$ generated by the cylinder sets. The second one is usually defined over the $\sigma$-algebra
of $\mathcal{B}_i= (S^1)^{\mathbb{L}}$ generated by the cylinder sets, where $\mathbb{L}$ is any countable set.
In order to show the relation of this two constructions in this paper, we focus on the cases where $\mathbb{L}=\mathbb{Z}$.

We will call $M_A$ the Gibbs-TF-$\mathbb{Z}$ for $A$, which is, by definition, the natural extension
of $m_A$, the Gibbs-TF for $A$.

For  a large class of potentials (see \cite{G}) we can show that $\mu^{ A,\sigma'}$ is independent of
the choice of $\sigma'\in  (S^1)^{\mathbb{N}}$. Here using a very simple argument we give a proof of this independence
using Ruelle operator when one consider free on the left, and a fixed $\sigma'\in (S^1)^{\mathbb{N}}$ boundary conditions.
We also show that this unique probability measure constructed using the Gibbs-DLR approach is equals to the measure $M_A$ obtained in the Gibbs-TF-$\mathbb{Z}$ for $A$. In a forthcoming paper we discus in great generality the equivalence of Gibbs-TF and Gibbs-DLR for one dimensional systems.

\subsection{Gibbs-DLR Probability Measures on $({\mathbb{S}^1})^\mathbb{Z}$.}

For a $\mathcal{B}_i$-measurable function $A:\mathcal{B}_i \to \mathbb{R}$ depending on the two first coordinates, we associated
a family $\Phi=(\Phi_{\Gamma})_{\Gamma\subset\mathbb{N}}$ of functions from $\mathcal{B}_i$ to $\mathbb{R}$, given by
$$
\Phi_{\Gamma}(x)=
\left\{
\begin{array}{rl}
 -A(x_n,x_{n+1}),& \text{if}\ \Gamma=\{n,n+1\};\\[0.3cm]
 0,&\text{otherwise}.
\end{array}
\right.
$$

We call this family $\Phi$ an interaction. For each $n\in\mathbb{N}$ we consider the associated Hamiltonian
\begin{equation}\label{def-Ham1}
H^\Phi_{\Lambda_n}(x)=-\sum_{k=-n}^{n-1}   A(x_{k},x_{k+1})\,,
\end{equation}
where $\Lambda_n=[-n,n]\cap\mathbb{Z}$.

The first step to obtain a Gibbs-DLR probability measure for a given $A:\mathcal{B}_i = ({\mathbb{S}^1})^\mathbb{Z}\to \mathbb{R}$ depending on the two first coordinates, whit boundary condition $\sigma' \in ({\mathbb{S}^1})^\mathbb{N}$,
is to construct a family of probability measures $
\mu_{\Lambda_n}^{\Phi,\sigma'}$ over $\mathcal{B}_i$ and then take
cluster points in the weak* topology of this family when $n\to\infty$.
Note that at least one cluster point exists because of the Banach-Alaoglu Theorem and any element on the set of these cluster points  will be called a Gibbs-DLR measure.
Once we take the limit when $n$ goes to infinity, the sequence of the sets $\Lambda_n= \{-n, -n+1,\ldots,-1,0,1,\ldots,n-1, n\}$ converges
in the set theoretical sense to $\mathbb{Z}$, which allows for these measures to capture information in the past and in the future coordinates.

Fixed a configuration $\sigma'=(\sigma'_0, \sigma'_1,.., \sigma'_n,..)\in(\mathbb{S}^1)^{\mathbb{N}}$, and, a potential $A$ as above, then, we define the Hamiltonian
on $\Lambda_n$ for the potential $\Phi$ with $\sigma'$ right boundary conditions by
$$H^\Phi_{\Lambda_n}(\tau|\sigma')=-\sum_{k=-n}^{n-2}   A(\tau_{k},\tau_{k+1})-A(\tau_{n-1},\sigma'_n)\,.$$

Note that $H^\Phi_{\Lambda_n}(\tau|\sigma')$ can also be considered as a function defined on $[0,2\pi]^{2n}$, i.e.,
$$H^\Phi_{\Lambda_n}(\tau_{-n},...,\tau_{n-1}|\sigma')=-\sum_{k=-n}^{n-2}   A(\tau_{k},\tau_{k+1})-A(\tau_{n-1},\sigma'_n)\,.$$

Let $M(\Lambda_n,\sigma')=\{ x \in ({\mathbb{S}^1})^\mathbb{Z}\,\, | \;x_i =  \sigma_i'\;,\; \forall \; i\geq n  \}$,
$d\nu$ the  Lebesgue probability measure on $\mathbb{S}^1$ (which we identify with $[0,2\pi]$)
and $d\nu_n$ is the Lebesgue probability measure on $(\mathbb{S}^1)^n$.

The partition function associated to the potential $\Phi$ with right boundary condition $\sigma' \in \mathcal{B}$
on the volume $\Lambda_n$ is defined by
\begin{eqnarray*}
Z_{\L_n}^{\Phi, \sigma'}&:=& \int_{M(\L_n,\sigma')}\,\,e^{- H_{\L_n}^{\Phi}(\tau|\sigma')}
\,d\nu(\tau)\\[0.2cm]
&=& \int_{[0,2\pi]^{2n}}\,\,e^{- H_{\L_n}^{\Phi}(\tau_{-n},...,\tau_{0}\,...,\tau_{n-1}|\sigma')}
\,d\nu(\tau_{-n})\,...\,d\nu(\tau_{0})\,...\, d\nu(\tau_{n-1}).
\end{eqnarray*}
We restrict our attention to potentials $\Phi$ for which the partition $Z_{\L_n}^{\Phi, \sigma'}$ is finite for
any choice $n$ and $\sigma'$. Hence for each $n$, this defines a probability measure which acts on continuous functions $f:\mathcal{B}\to\mathbb{R}$ (depending on finite coordinates) by
$$
\int_{\mathcal{B}_i} f \,\,d\mu_{\Lambda_n}^{\Phi,\sigma'}=\, \frac{1}{ Z_{\L_n}^{\Phi,  \sigma'}}
\,
\int_{M(\Lambda_n,\sigma')}\!\!\!\!\!\!f(\tau)\, e^{- H_{\Lambda_n}^{\Phi}(\tau|\sigma')} d\nu(\tau_{-n})d\nu(\tau_{-n+1})\,...\, d\nu(\tau_{n-1})\,.
$$

Note that in this way for any fixed $\sigma'$ the probability measure $\mu_{\Lambda_n}^{\Phi,\sigma'}$ depends just on $A$ (and, of course, $\sigma'$), thus we could also denoted it $\mu_{\Lambda_n}^{A,\sigma'}$. But here we will adopt the Statistical Mechanics notation $\mu_{\Lambda_n}^{\Phi,\sigma'}$ as used in \cite{G} and \cite{Si}.

For a fixed $\sigma'$ we are interested in the limit of $\mu_{\Lambda_n}^{\Phi,\sigma'}$, when $n \to \infty$.
Any possible cluster point of this sequence will be denoted by $\mu^{A,\sigma'}$ (or, $\mu^{\Phi,\sigma'}$).
Any one of these is called  a Gibbs state for $A$  with a boundary condition $\sigma' \in\mathcal{B} $ on the right and free on the left.

Given $A:\mathcal{B} \to \mathbb{R}$,
by the major theorem  of section \ref{sec1}, we know there is a maximal positive eigenvalue $\lambda=\lambda_A$ associated to the eigen-function $\psi_A$.
We also have, for any $\psi: \mathcal{B} \to \mathbb{R}$,
\begin{equation}\label{ruelleiterado}
\mathcal{L}_{ A}^n\psi(y)=
\int_{[0,2\pi]^n} e^{ S_n A(\tau y)}\psi(\tau y) \ d\nu_n(\tau) \,.
\end{equation}
If $A$ depends on two coordinates, then, $\psi_A$ depends on one coordinate (as  we get from section \ref{sec3}).
Note that for any $\tau \in ({\mathbb{S}^1})^\mathbb{N}$
we have
$\mathcal{L}_{A}^{2 n} (1)\,(\sigma^{n} (\tau)) =
Z_{\L_n}^{\Phi,\tau},$
where $\sigma$ is the shift on $({\mathbb{S}^1})^\mathbb{N}$ and $\mathcal{L}_{\bar A}^n 1=1$ for any $n \in \mathbb{N}$, where
$$\bar A = A + \log \psi_A - \log \psi_A \circ \sigma - \log \lambda_A.$$

Let $\bar \Phi$ be the potential defined by $\bar A$ and $\pi$ the natural  projection of  $({\mathbb{S}^1})^\mathbb{Z}$ to  $({\mathbb{S}^1})^\mathbb{N}$. (analogous to the case for the potential $A$), we set  for any Borelian $C \subset \mathcal{B}$
$$
\pi\, \mu_{\Lambda_n}^{\bar \Phi,\sigma'}(C)=
\frac{1}{Z_{\L_n}^{\bar \Phi,\sigma'}}
\,
\int_{M(\L_n,\sigma')}\!\!\!\!\!\!\mathds{1}_C(\tau) e^{- H_{\L_n}^{\bar \Phi}(\tau|\sigma')} d\nu(\tau)\,.
$$

We point out that a potential $A$ which depends on two coordinates can be seen as a potential defined either in  $({\mathbb{S}^1})^\mathbb{N}$, or  $({\mathbb{S}^1})^\mathbb{Z}$. Another important remark is when the spin variables take values in the close interval $[-1,1]$ these models are known in the literature as continuous Ising model.

\begin{proposition} Consider a fixed $\sigma'\in \mathcal{B}= ({\mathbb{S}^1})^\mathbb{N}$
Given $A:\mathcal{B}_i\to \mathbb{R} $, which depends on two coordinates, if $\bar A$
is its normalized associated potential then for any cluster point $\pi \mu^{ \bar\Phi,\sigma'}$ we have that
$$
m=  \pi \mu^{ \bar\Phi,\sigma'},
$$
where $m=m_{\bar A}$ is the Gibbs-TF measure for $\bar A$.
\end{proposition}
We will show that  $\displaystyle\lim_{n \to \infty} \pi \mu_{\Lambda_n}^{\bar \Phi,\sigma'}=m$, so this limit does not  depend on the fixed $\sigma'$ we choose.

\noindent{\bf Proof.} Consider a given  $f: \mathcal{B} \to \mathbb{R} $
which depends on finitely many coordinates, (let's say $r>0$).
Note that
$$H^{\bar\Phi}_{\Lambda_n}(\tau|\sigma')= - \sum_{k=-n}^{n-2}   \bar A (\tau_{k},\tau_{k+1})-\bar A(\tau_{n-1},\sigma'_n)\,,$$
and that $Z_{\Lambda_n}^{\bar\Phi,\sigma'}=1$. Suppose that $n>r$. By definition
\begin{eqnarray*}
\int f \,\,d\, \pi\, \mu_{\Lambda_n}^{\bar\Phi,\sigma'}&=&
\,
\int_{M(\Lambda_n,\sigma')}\!\!\!\!\!\!f(\tau)\, e^{- H_{\Lambda_n}^{\bar\Phi}(\tau|\sigma')} d\nu(\tau)
\\[0.3cm]
&=& \int_{[0,2\pi]^{2n}}\!\!\!\!\!\!f(\tau_0,...,\tau_r)\, e^{- H_{\Lambda_n}^{\bar\Phi}(\tau_{-n}...\tau_{n-1}|\sigma')} d\nu(\tau_{-n})...d\nu(\tau_{n-1})\,
\\[0.3cm]
&=& \int_{ [0,2\pi]^{n}}f(\tau_0,...,\tau_r)
e^{\sum_{k=0}^{n-2}\bar A(\tau_k,\tau_{k+1}) + \bar A(\tau_{n-1},\sigma_n')} \times
\\
&&\times
\left(\int_{ [0,2\pi]^{n}}\, e^{\sum_{k=-n}^{-1}\bar A(\tau_k,\tau_{k+1})} d\nu(\tau_{-n})...d\nu(\tau_{-1})\right)d\nu(\tau_{0})\,...\, d\nu(\tau_{n-1})
\\[0.3cm]
&=& \int_{ [0,2\pi]^{n}}f(\tau_0,...,\tau_r)
e^{\sum_{k=0}^{n-2}\bar A(\tau_k,\tau_{k+1}) + \bar A(\tau_{n-1},\sigma_n')}d\nu(\tau_{0})\,...\, d\nu(\tau_{n-1}).
\end{eqnarray*}
where in the last equation we used $n$ times $\int_{[0,2\pi]} e^{\bar A(x,y)} d\nu(x) =1$.


In this way,

$$
\int f \,\,d\, \pi\, \mu_{\Lambda_n}^{\bar\Phi,\sigma'}
=\mathcal{L}_{\bar A}^{ n} (f) (\sigma^n  (\sigma ')).$$

Is is known from section \ref{sec1} that $\mathcal{L}_{\bar A}^{\,n} (f)$ converges uniformly to $ \int f \, dm$, as $n$ goes to
infinity, where $m$ is Gibbs-TF for $A$ or $(\bar A$).
As the convergence of $\mathcal{L}_{\bar A}^n (f)$, when $n \to \infty$, is uniform, then
$\lim_{n \to \infty} \pi \mu_{\Lambda_n}^{\bar\Phi,\sigma'}=m$

\cqd

\bigskip

\begin{corollary}
 For any $\sigma' \in  ({\mathbb{S}^1})^\mathbb{N}$, and, any $f$ which depends on finitely many coordinates

$$\frac{
\int_{ M(\Lambda_n,\sigma')}f(\tau)\, e^{- H_{\Lambda_n}^{\bar\Phi}(\tau|\sigma')} d\nu(\tau_{-n})d\nu(\tau_{-n+1})\,...\, d\nu(\tau_{n-1})}{ \int_{(\mathbb{S}^1)^{\L_n}}f(\tau)\, e^{- H_{\Lambda_n}^{\bar\Phi}(\tau)}
d\nu(\tau_{-n})d\nu(\tau_{-n+1})\,...\, d\nu(\tau_{n-1} ) d\nu(\tau_{n}) } \to 1 \,,
$$
when $n \to \infty$.
\end{corollary}

\noindent{\bf Proof.} This follows easily from the above because the convergence of $\mathcal{L}_{\bar A}^{\,n} (f)$ is uniform.

\cqd

\begin{proposition}Suppose $\sigma'\in ({\mathbb{S}^1})^\mathbb{N}$.
Given $A:\mathcal{B}_i\to \mathbb{R} $, which depends on two coordinates
and, a coboundary $h:\mathcal{B}_i\to \mathbb{R}$, which depends on one  coordinate (the $0$ coordinate), and,  such that
$$\bar A = A +h - h \circ \hat{\sigma} + \log \lambda,$$
where $\hat{\sigma} $ is the shift on $\mathcal{B}_i$, then
$$
\pi (\mu^{\Phi,\sigma'})=  \pi (\mu^{\bar\Phi,\sigma'}).
$$
\end{proposition}
\noindent{\bf Proof.} Consider a function $f: \mathcal{B}\to \mathbb{R}$ which depends on finite coordinates $f(\tau_0,\tau_1,..,\tau_k)$, $k>0$. We have first that
\begin{eqnarray*}
H^{\bar\Phi}_{\Lambda_n}(\tau|\sigma')&=&
-\bar A ( \tau_{-n} ,\tau_{-n+1}) -\sum_{k=-n+1}^{n-2}   \bar A (\tau_{k},\tau_{k+1})-\bar A(\tau_{n-1},\sigma'_n)\,\\[0.3cm]
&=& H^{\Phi}_{\Lambda_n}(\tau|\sigma')+ h(\sigma'_n) - h(\tau_{-n})-2\, n\,\log \lambda.
\end{eqnarray*}
Hence
$$
-H^{\Phi}_{\Lambda_n}(\tau|\sigma')=
-H^{\bar\Phi}_{\Lambda_n}(\tau|\sigma')+ h(\sigma'_n) - h(\tau_{-n})-2\, n\,\log \lambda.
$$
Therefore
$$\int_{M(\Lambda_n,\sigma')}\,f(\tau)\, e^{- H_{\Lambda_n}^{\Phi}(\tau|\sigma')}\,d\nu(\tau)=$$
 $$\lambda^{-2 n} \, \, e^{ h(\sigma_n ') }\int_{M(\Lambda_n,\sigma')}\,e^{-h(\tau_{-n})}\, f(\tau_0,..,\tau_k)\,
 e^{- H_{\Lambda_n}^{\bar\Phi}(\tau|\sigma')}\,d\nu(\tau_{-n})\,...\, d\nu(\tau_{n-1}),$$
by taking $f=1$ we have
$$Z_{\L_n}^{\Phi,\sigma'}=\int_{M(\Lambda_n,\sigma')}\, e^{- H_{\Lambda_n}^{\Phi}(\tau|\sigma')}\,d\nu(\tau)=$$
 $$\lambda^{-2 n} \, \, e^{ h(\sigma_n ') }\int_{M(\Lambda_n,\sigma')}\,e^{-h(\tau_{-n})}\, \, e^{- H_{\Lambda_n}^{\bar\Phi}(\tau|\sigma')}\,d\nu(\tau_{-n})\,...\, d\nu(\tau_{n-1}).$$
$$ 
=\lambda^{-2 n} \, \, e^{ h(\sigma_n ') }\mathcal{L}_{\bar A}^{2 n} (e^{-h})( \sigma^n (\sigma')).$$

We already shown in the previous sections that

$$
\mathcal{L}_{\bar A}^{2 n} (e^{-h})( \sigma^n (\sigma'))\, \to \int\, e^{-h}\,\, dm_{\bar A},
$$
uniformly in $n$.
Therefore, $Z_{\L_n}^{\,\Phi,\sigma'} \sim\, \lambda^{-2 n}\, e^{ h(\sigma_n ')}   \int\, e^{-h}\,\, dm_{\bar A}$.


 We also have
$$\int_{M(\Lambda_n,\sigma')}\,e^{-h(\tau_{-n})}\, f(\tau_0, \tau_1,..,\tau_k)\, e^{- H_{\Lambda_n}^{\bar\Phi}(\tau|\sigma')}\,d\nu(\tau)=$$

 $$ \int_{ [0,2\pi]^{n}}f(\tau)
e^{\sum_{k=0}^{n-2}\bar A(\tau_k,\tau_{k+1}) + \bar A(\tau_{n-1},\sigma_n')}
\,\,\times$$
$$\times\left(\int_{ [0,2\pi]^{n}}\,e^{-h(\tau_{-n})} e^{\sum_{k=-n}^{-1}\bar A(\tau_k,\tau_{k+1})} \prod_{i=-n}^{-1}d\nu(\tau_i)\right)\prod_{k=0}^{n-1}d\nu(\tau_k)=  $$
$$ \int_{ [0,2\pi]^{n}}f(\tau)
e^{\sum_{k=0}^{n-2}\bar A(\tau_k,\tau_{k+1}) + \bar A(\tau_{n-1},\sigma_n')}
(\mathcal{L}_{\bar A}^{ n} (e^{-h})( \tau))\prod_{k=0}^{n-1}d\nu(\tau_k)= $$

 $$  \int_{ [0,2\pi]^{n}}f(\tau)
e^{\sum_{k=0}^{n-2}\bar A(\tau_k,\tau_{k+1}) + \bar A(\tau_{n-1},\sigma_n')}
\bigg(\mathcal{L}_{\bar A}^{ n} (e^{-h})( \tau)-\int e^{-h} dm_{\bar A}\bigg)\prod_{k=0}^{n-1}d\nu(\tau_k)+$$
 $$+\int_{ [0,2\pi]^{n}}f(\tau)
e^{\sum_{k=0}^{n-2}\bar A(\tau_k,\tau_{k+1}) + \bar A(\tau_{n-1},\sigma_n')}
\left(\int e^{-h} dm_{\bar A}\right)\prod_{k=0}^{n-1}d\nu(\tau_k)  $$
$$\to \int\, f d m_{\bar A} \int e^{-h} dm_{\bar A}.$$
where in the convergence we used the fact that,
given any $\epsilon>0$, there exists $N_{\epsilon}$ such that, for $n>N_{\epsilon}$, we have

 $$  \left|\int_{ [0,2\pi]^{n}}f(\tau)
e^{\sum_{k=0}^{n-2}\bar A(\tau_k,\tau_{k+1}) + \bar A(\tau_{n-1},\sigma_n')}
\bigg(\mathcal{L}_{\bar A}^{ n} (e^{-h})( \tau)-\int e^{-h} dm_{\bar A}\bigg)\prod_{k=0}^{n-1}d\nu(\tau_k)\right| < $$
 $$ <\epsilon \int_{ [0,2\pi]^{n}}f(\tau)
e^{\sum_{k=0}^{n-2}\bar A(\tau_k,\tau_{k+1}) + \bar A(\tau_{n-1},\sigma_n')}
\prod_{k=0}^{n-1}d\nu(\tau_k)=
$$
$$ = \epsilon \mathcal{L}_{\bar A}^{ n} (f)( \sigma^n (\sigma'))<2 \epsilon \int fdm_{\bar A}\,,
$$
which means that the first integral vanishes when $n \to \infty$, while the second integral is
$$\int_{ [0,2\pi]^{n}}f(\tau)
e^{\sum_{k=0}^{n-2}\bar A(\tau_k,\tau_{k+1}) + \bar A(\tau_{n-1},\sigma_n')}
\left(\int e^{-h} dm_{\bar A}\right)
\prod_{k=0}^{n-1}d\nu(\tau_k)
=$$
$$= \mathcal{L}_{\bar A}^{ n} (f)( \sigma^n (\sigma')) \int e^{-h} dm_{\bar A}\to \int\, f d m_{\bar A} \int e^{-h} dm_{\bar A}
 \,.$$




 Finally,

 $$\frac{\int_{M(\Lambda_n,\sigma')}\,f(\tau)\, e^{- H_{\Lambda_n}^{\Phi}(\tau|\sigma')}\,d\nu(\tau)}
 {Z_{\L_n}^{\Phi,\sigma'}}
 $$
 $$=
 \frac{\lambda^{-2 n} \, \, e^{ h(\sigma_n ')}\int_{M(\Lambda_n,\sigma')}e^{-h(\tau_{-n})}\,f(\tau)\,
 e^{- H_{\Lambda_n}^{\bar\Phi}(\tau|\sigma')}\,d\nu(\tau)}
 {Z_{\L_n}^{\Phi,\sigma'}}
 $$
 $$\sim\frac{\int_{M(\Lambda_n,\sigma')}e^{-h(\tau_{-n})}\,f(\tau)\, e^{- H_{\Lambda_n}^{\bar\Phi}(\tau|\sigma')}\,d\nu(\tau)}{\int e^{-h} dm_{\bar A}}\to\int\, f d m_{\bar A}$$

Therefore, $\pi \mu^{\Phi,\sigma'}=  \pi \mu^{\bar\Phi,\sigma'}$.

\cqd

\begin{corollary}
Consider a general  $\sigma'\in \mathcal{B}$.
Given $A:\mathcal{B}\to \mathbb{R} $, then,
$$m_A=  \pi \mu^{\Phi,\sigma'},$$
where $m=m_{A}$ is the Gibbs-TF for $ A$.
\end{corollary}

\noindent{\bf Proof.} It follows from $\bar A = A + \log \psi_A - \log \psi_A \circ \sigma - \log \lambda_A.$

\cqd

\bigskip

According to \cite{G} Part III page 289, for any $\sigma'$
the probability measure $\mu^{  A,\Phi,\sigma'}$ is invariant for $\hat{\sigma}$ acting on  $({\mathbb{S}^1})^\mathbb{Z}$.

By definition, the Gibbs-FT-Z state $M_A$ on $({\mathbb{S}^1})^\mathbb{Z}$, is the natural extension of $m_A$, and, it is also invariant for $\hat{\sigma}$ acting on  $({\mathbb{S}^1})^\mathbb{Z}$.

\bigskip

\begin{proposition} Suppose $A: ({\mathbb{S}^1})^\mathbb{Z}\to \mathbb{R}$ depends on two coordinates, and, consider $\sigma'\in \mathcal{B}$,
then
$$\mu^{\Phi,\sigma'}= M_{ A}.$$
\end{proposition}

\noindent{\bf Proof.}
$\mu^{\Phi,\sigma'}$ and $M_{A}$ are both the natural extension of $m_A$.

\cqd
\bigskip

\begin{proposition} Suppose $A: ({\mathbb{S}^1})^\mathbb{Z}\to \mathbb{R}$ depends on two coordinates, and, consider $\sigma',\sigma''\in \mathcal{B}$,
then
$$\mu^{\Phi,\sigma'}= \mu^{\Phi,\sigma''}.$$
\end{proposition}

\noindent{\bf Proof.}
$\mu^{\Phi,\sigma'}$ and $  \mu^{\Phi,\sigma''}$ are both the natural extension of $m_A$.

\cqd
\bigskip

The final conclusion is that, if the potential depends on two coordinates, then the Gibbs probability measure on $({\mathbb{S}^1})^\mathbb{Z}$ in both settings, Thermodynamic Formalism and Statistical Mechanics via a boundary condition $\sigma'$ on the right side, coincide.

\bigskip

Now we will analyze the free-boundary case. Remember that
$$H^\Phi_{\Lambda_n}(\tau)=-\sum_{k=-n}^{n-1}   A(\tau_{k},\tau_{k+1})\,.$$

We are going to define the Gibbs probability measure in the sense of Statistical
Mechanics with free boundary condition on the left and on the right. For a given $n>0$,
$$
Z_{\L_n}^{\Phi}= \int_{(\mathbb{S}^1)^{\L_n}}\,\,
e^{- H_{\L_n}^{\Phi}(\tau)} \,d\nu(\tau_{-n})d\nu(\tau_{-n+1})\,...\, d\nu(\tau_{n-1}) d\nu(\tau_{n})
$$
will be the partition function which
corresponds to the case of free a boundary condition  on the right and on the left.

For each $n$, this defines a probability measure which acts on continuous functions $f$ (depending on finite coordinates) by
$$
\int f \,\,d\, \, \mu_{\Lambda_n}^{\Phi}=\, \frac{1}{ Z_{\L_n}^{\Phi}}
\,
\int_{(\mathbb{S}^1)^{\L_n}}\!\!\!\!\!\!f(\tau)\, e^{- H_{ A,\Lambda_n}^{\Phi}(\tau)} d\nu(\tau_{-n})d\nu(\tau_{-n+1})\,...\, d\nu(\tau_{n-1})\,d\nu(\tau_{n}).
$$

Any weak limit  of subsequences of  $\mu_{\Lambda_n}^{\Phi}$ will be called a
Gibbs state for $A$ with  a free  boundary condition  on the right and on the left.

It follows from Corollary 1 above that any Gibbs state for $A$ with  a free  boundary condition  on the right and on the left is equal to $M_A$.

The result we will analyze in the next section will be the case of  a free  boundary condition  on the right and on the left.

\vspace{0.6cm}

\section{An example by A. C. D. van Enter and W. M. Ruszel where there is no selection}\label{sec5}

In this section we will consider $A$ depending on its first neighbors, and having the form $A(x)=A(x_0,x_1)= U(x_0-x_1)$.

We want to show  a particular example (introduced by \cite{van}), where the potential is not continuous and is of the form: $\tilde U:[0,2\pi] \rar \R$ is a function such that $\tilde U|_{[a_n,b_n)}$, is constant for each $n$ and equal to $c_n$, where $[a_n,b_n)$, $n \in \mathbb{N}$ is a partition of $[a,b]$.

We will show that  for each positive $\beta$ we can  also consider an extension of Gibbs-TF, say $\mu_{\beta,\tilde U}$, over ${\cal B}$ and also
that this measure coincides with the Gibbs-DLR for this potential $\tilde U$.
 In \cite{van} the authors have shown that there is no selection of the family $\mu_{\beta, \tilde U}$ when $\beta \to \infty$.

We will present here all the details of the proof of this non-trivial result.

Basically, we will show that $\int I_B \,d \mu_{\beta,\tilde U}$ does not converge when $\beta \to \infty$, for a set $B$
which depends just on the coordinates $(x_0,x_1)$.
Therefore, this is also the same as to say that  $\int I_B\, d \hat{\mu}_{\beta,\tilde U}$
does not converge (see Remark 4 just before proposition \ref{proponly}).

The main result of this section is theorem \ref{nonselec-theo}, which is a consequence of  corollary \ref{nonselec-cor}
and  lemma \ref{nonselec-lemma}. Subsection \ref{VanEntercomocasolimite} shows that results of previous sections are still valid even if the potential $A$ belongs to certain classes of non-continuous potentials including the potential of \cite{van}.

\bigskip



\subsection{Gibbs Measures for Non-continuous Potentials and DLR formulation of Statistical Mechanics}\label{VanEntercomocasolimite}

So far we have defined Gibbs Measures for Holder continuous potentials in sections \ref{sec1} (general case) and \ref{sec3} (nearest neighbors interaction, i.e. potential depending on two coordinates). In the section \ref{sec3} we gave an alternative definition based on transition kernels associated to a certain potential (or Hamiltonian) $A$, and proved that this definition is equivalent to the one of section \ref{sec1}.

We will now show that our definition coincides with the usual one in
Statistical Mechanics,
in the case of a certain special
non-continuous potential depending on two coordinates. We assume, among other things, the form $\tilde A(x) =\tilde A(x_0,x_1)=\tilde U(x_0-x_1)$, where $\tilde U:\mathbb{S}^1 \to \mathbb{R}$ is a bounded $L^1$ function, which is pointwise approximated by Holder functions  $U_n$. This case will cover the important example to be described later. A potential of this form is called symmetric.


First we will show that the main results of Section \ref{sec3} are true for this potential $\tilde A(x) =\tilde U(x_0-x_1)$, which is no longer continuous.

Using the notation described in section \ref{sec3}, let
${\cal L}_{\beta \tilde U},{\bar{\cal L}}_{\beta \tilde U}:C([0,2\pi])\to C([0,2\pi])$ be given by

\begin{equation}\label{LL}   {\cal L}_{\beta \tilde U} \psi (y) = \frac{1}{2\pi}\int_0^{2\pi} e^{\beta \tilde U(x-y)} \, \psi(x)dx,  \end{equation}

\begin{equation}\label{LLadj}   {\bar{\cal L}}_{\beta \tilde U} \psi (x) =\frac{1}{2\pi} \int_0^{2\pi} e^{\beta \tilde U(x-y)} \, \psi(y)dy.  \end{equation}
for any $y \in [0, 2 \pi]$.

In order to simplify the notation we denote $  {\cal L}_{\beta}$ instead of  $ {\cal L}_{\beta \tilde U}$ .

\begin{lemma}\label{preservacone}
The operators ${\cal L}_{\beta}$ and ${\bar{\cal L}}_{\beta}$ preserve the set of continuous functions in $[0,2\pi]$, sending continuous functions to uniformly continuous functions. Moreover, a bounded function is mapped to an uniformly continuous one.
\end{lemma}

The fact  that continuous functions are preserved implies the compactness of the operator, as we can see in pages 43 and 47 of \cite{Co}.

\proof Consider a fixed $\beta$ and the operator ${\cal L}_{\beta}$.
Let $f$ be a continuous function.

Fix $\epsilon >0$. Let $A_c$ be a continuous function such that $$\|\tilde A-A_c\|_{L^1}< \frac{\epsilon}{4 \|f\|_{C^0}}\,.$$ Here we use the $L_1$ norm on the functions defined on the one-dimensional set $[0, 2\pi]$.
 Such a function exists because continuous functions are dense in $L_p[0,2\pi]$ for $p\geq 1$.

Let $K_c(x,y)=A_c(x-y)$. 
We have $\tilde A=A_c+(\tilde A-A_c)$

Moreover, let $\delta>0$ be such that $$|A_c(z)-A_c(w)|< \frac{\epsilon}{2 \|f\|_{C^0}}$$ if $|z-w|<\delta$.

Suppose $|y_1-y_2|< \delta$. Then we have
$$ |\LL(f)(y_1) - \LL(f)(y_2)| =  \left|\int K(x,y_1) f(x) dx - \int K(x,y_2) f(x) dx\right| $$
$$ \leq \int  \left|A_c(x-y_1) -  A_c(x-y_2)\right| |f(x)| dx + $$
$$ +\int |(\tilde A-A_c)(x-y_1)| |f(x)| dx
+ \int |(\tilde A-A_c)(x-y_2)| |f(x)| dx < \epsilon
\,.$$
\cqd

\bigskip

The proof of the next theorem is a small modification of the proof of Theorem 3 of \cite{LMST}

\begin{theorem}\label{KreinRuthman} 
The operators  ${\cal L}_{\beta}$ and ${\bar{\cal L}}_{\beta}$ have the same  positive maximal eigenvalue $\lambda_{\beta}$, which is simple and isolated. The eigenfunctions  associated to these operators, say $\psi_{\beta}$ and $\bar\psi_{\beta}$, are positive functions.
\end{theorem}

\noindent{\bf Proof.}
We can see that
${\cal L}_{\beta}$ is a compact operator, because Lemma \ref{preservacone} shows that the image of the unity closed ball of $C([0,1])$ under ${\cal L}_{\beta}$ is an
equicontinuous family in $C([0,1])$. Thus, we can use Arzelà-Ascoli theorem to prove the compactness of ${\cal L}_{\beta}$ (see also  Chapter IV, section 1 of \cite{Sch}).

The spectrum of a compact operator contains a sequence of eigenvalues that converges to zero, possibly contains zero. This implies that any non-zero eigenvalue of ${\cal L}_{\beta}$ is isolated
(i.e. there is no sequence in the spectrum of ${\cal L}_{\beta}$ which converges to a non-zero eigenvalue).

The definition of ${\cal L}_{\beta}$ now shows that ${\cal L}_{\beta}$ preserves the cone of positive functions in $C([0,1])$, sending a point in this cone to the interior of the cone. This means that ${\cal L}_{\beta}$ is a  positive operator.

The Krein-Ruthman theorem (Theorem 19.3 of \cite{De}) implies that there exists a
 positive eigenvalue $\lambda_{\beta}$, which is  maximal (i.e. if $\lambda\neq\lambda_{\beta}$ is in the spectrum of ${\cal L}_{\beta}$ then $\lambda_{\beta}>|\lambda|$.)
and simple  (i.e. the eigenspace associated to $\lambda_{\beta}$ is one-dimensional). Moreover $\lambda_{\beta}$ is associated to a positive eigenfunction
$\psi_{\beta}$. 

If we proceed in the same way as in \cite{LMST}, we obtain the same conclusions about the operator ${\bar{\cal L}}_{\beta}$, and we get
the respective eigenvalue $\bar\lambda_{\beta}$ and eigenfunction $\bar{\psi}_{\beta}$.

In order to prove that $\bar\lambda_{\beta}= \lambda_{\beta}$, we use
the positivity of $\psi_{\beta}$ and $ \bar{\psi}_{\beta}$ and the fact that ${\bar{\cal L}}_{\beta}$
is the adjoint of ${{\cal L}}_{\beta}$. (Here we see that  our operators can be, in fact, defined in the Hilbert space $L^2([0,1])$, which contains $C([0,1])$ ). We have
$<\psi_{\beta},\bar{\psi}_{\beta}> = \int \psi_{\beta}(x) \bar{\psi}_{\beta}(x) dx > 0$, and
$$
\lambda_{\beta} <\psi_{\beta},\bar{\psi}_{\beta}> = <{{\cal L}}_{\beta}\psi_{\beta},\bar{\psi}_{\beta}> =
<\psi_{\beta},{\bar{\cal L}}_{\beta}\bar{\psi}_{\beta}> = \bar\lambda_{\beta} <\psi_{\beta},\bar{\psi}_{\beta}>.$$ \cqd

 By the periodicity of $\tilde U$, ${\cal L}_{\beta \tilde U} \psi (1) $ and $ {\bar{\cal L}}_{\beta \tilde U}(1) $ are independent of $x$. Therefore $\psi_{\beta,\tilde U}(x)=\bar\psi_{\beta,\tilde U}(x)=1$ are the eigenfunctions associated to the maximal eigenvalue $\lambda_{\beta,\tilde U}$.

It is easy to see that  $$\lambda_{\beta,\tilde U}=\frac{1}{2\pi} \int_0^{2\pi} e^{\beta \tilde U(x-y)} \, dy\,.$$
In the notation section \ref{sec3},  $\theta_{\beta,\tilde U}(x)=1$ and the transition Kernel is given by

$$K_{\beta, \tilde U}(x,y):=\frac{e^{\beta \tilde U(x-y)}\,\,}{\,\lambda_{\beta}}\,. $$

For instance, for any cylinder
$$\mu_{\beta,\tilde U}( A_0...A_k)=\int_{A_0...A_k}
\, \frac{e^{\beta\sum_{i=0}^{k-1} \tilde U(x_i-x_{i+1})}\,\,}{\,\lambda_{\beta}^k}\;dx_k\,... \,dx_0.$$

This measure does not came from a H\"older potential, but we can approximate this measure by Gibbs-TF measures associated to H\"olders potentials, as we will see next.

Let us now analyze the case $ A(x) = A(x_0,x_1)= U(x_0-x_1)$ where $U:\R \rar \R$ is a Holder continuous function  $2\pi$-periodic.
By the same arguments used above, it is easy to see that $\psi_{\beta,U}(x)=\bar\psi_{\beta,U}(x)=1$ are the the eigenfunctions of the operators ${\cal L}_{\beta, U},{\bar{\cal L}}_{\beta ,U}$ associated to the maximal eigenvalue $\lambda_{\beta,U}$ (see section \ref{sec3}), where

$$\lambda_{\beta,U}=\frac{1}{2\pi} \int_0^{2\pi} e^{\beta U(x-y)} \, dy\,.$$

As in section \ref{sec3},  $\theta_{\beta,U}(x)=1$ and the transition Kernel is given by

$$K_{\beta, U}(x,y):=\frac{e^{\beta U(x-y)}\,\,}{\,\lambda_{\beta}}\,. $$

Hence, for any cylinder
$$\mu_{\beta,U}( A_0...A_k)=\int_{A_0...A_k}
\, \frac{e^{\beta\sum_{i=0}^{k-1} U(x_i-x_{i+1})}\,\,}{\,\lambda_{\beta}^k}\;dx_k\,... \,dx_0.$$
By theorem \ref{spitzer_vs_PP} we see  that  $\mu_{\beta,U}=m_U$, the Gibbs-TF for $U$.

Let now $\tilde U$ be a $L^1$ potential such that there exists an uniformly bounded  sequence of Holder continuous potentials $U_n$ converging point wise to $\tilde U$.

By the Dominated Convergence Theorem, we have that
$$\lambda_{\beta,U_n}=\frac{1}{2\pi} \int_0^{2\pi} e^{\beta U_n(x-y)} \, dy \to \frac{1}{2\pi} \int_0^{2\pi} e^{\beta \tilde U(x-y)} = \lambda_{\beta,\tilde U}
\,,$$
as $k \to \infty$,
and also for any cylinder $A_0...A_k$, we have

$$\mu_{\beta,U_n}( A_0...A_k)  \to \int_{A_0...A_k}
\, \frac{e^{\beta\sum_{i=0}^{k-1} \tilde U(x_i-x_{i+1})}\,\,}{\,\lambda_{\beta}^k}\;dx_k\,... \,dx_0=\mu_{\beta,\tilde U}( A_0...A_k)
$$
as $k \to \infty$.




Note that the measure $\mu_{\beta,\tilde U}$ coincides with the Gibbs-DLR measure of statistical mechanics in the special case of nearest neighbors interaction of the kind
 $\tilde A(x) =\tilde A(x_0,x_1)=\tilde U(x_0-x_1)$ as can be seen, for example, in
  \cite{G}. 
We also remark that $\frac{1}{\lambda_{\beta}^k}$ is the partition function of DLR formulations of statistical mechanics.

\vspace{.6cm}


\subsection{ One-dimensional Systems With Symmetric Potentials}\label{sec52}

To explain the no selection measure theorem we will use the formalism introduced in last section. Here we take $\mathbb{L}=\mathbb{Z}$ which is the origin of the term ``one-dimensional" in the title of this section. We assign for each $i\in \mathbb{Z}$ the measure space $(\mathbb{S}^1,\mathscr{B},\nu)$, where $\nu$ is the  Lebesgue probability  measure on the circle. For each $n\in\N$ we denote by $\Lambda_n=:[-n,n]\cap\Z$.
We will use free conditions on the left and on the right side.

For convenience, we use the natural measure isomorphism between
the Bernoulli spaces $(\mathbb{S}^1)^{\Z}$ and $[0,2\pi)^{\Z}$ to define the Hamiltonian we introduced before
in $(\mathbb{S}^1)^{\Z}$.
Let  $\Phi=(\Phi_{\Gamma})_{\Gamma\subset \mathcal{L}}$ be a family of functions on $[0,2\pi)^{\Z}$, such that
$$
\Phi_{\Gamma}(\theta)=
\left\{
\begin{array}{ll}
U(\theta_k-\theta_{k+1}),&\text{if}\ \ \Gamma=\{k,k+1\};\\[0.2cm]
0,&\text{otherwise.}
\end{array}
\right.
$$
where $U$ is a potential defined by the $2\pi$-periodic extension of
$$
\tilde{U}(x)=\sum_{j=1}^{\infty}c_j\mathds{1}_{A_{j}}(x)\,,
$$
and $\{A_j\}_{j\geq 1}$ is a partition of
$[0,2\pi)$ given by intervals of the form $A_j=[a_j,b_j)$.


Using the isomorphism and the family $\Phi$ mentioned above, the Hamiltonian  in the finite volume $\Lambda_n$, with boundary condition $x'$, will be given by the following expression, if $x=(\theta_{k})_{k \in \mathbb{Z}}$
\begin{equation}
\label{Hamiltonian}
H_{\Lambda_n}(x_{\Lambda_n}x'_{\Lambda_n^c})
=-\sum_{k=-n}^{n-1}U(\theta_k-\theta_{k+1})
-U(\theta_{-n}-\theta'_{-n-1})-U(\theta_n-\theta'_{n+1}).
\end{equation}

The family $\Phi$ we are considering is associated to a potential $A$ which depends only on the nearest-neighbors and given by $A(x,y)=U(x-y)$.
We can prove (see \cite{Ruelle}) that for each fixed $\beta\in(0,+\infty)$, the set
$\mathcal{G}_{\beta,\Phi}$ is a singleton set and its unique measure denoted by $\mu_{\beta A}$ is given by
$$
\mu_{\beta A}=
w-\lim_{\L_n\nearrow\N} \mu_{\Lambda_n}^{\beta,\Phi},
$$
where for all $n\in\N$ and $E\in\mathcal{B}_i$
\begin{equation}\label{med-gibbs-van-enter-vol-n}
\mu_{\Lambda_n}^{\beta,\Phi}(E)=
\frac{1}{Z_{\Lambda_n}^{\beta,\Phi}}
\int_{(\mathbb{S}^1)^{\Lambda_n}}
\!\!\!\!\!\!\mathds{1}_{\pi_{\Lambda_n}(E)}(\theta)
\exp\left(\!\!\beta\sum_{k=-n}^{n-1}U(\theta_k-\theta_{k+1})\!\!\right)
\!\!\! \prod_{k=-n}^{n}\!\!\! d\nu(\theta_k).
\end{equation}

From now on we call $\mu_{\Lambda_n}^{\beta,\Phi}$ the Gibbs measure in the volume $\Lambda_n$ for the Hamiltonian (\ref{Hamiltonian}) at inverse temperature $\beta$.
\bigskip

We are using above free boundary conditions on the left and on the right side.
\bigskip

We will consider here a real parameter $\beta$, wich means the inverse of the temperature, and the Gibbs probability measure $\hat{\mu}_{\beta A}$ over
$ {\cal B}_i=(\mathbb{S}^1)^\mathbb{Z}$ (see the considerations just above proposition \ref{proponly}).

Note that if $U$ has a unique maximum at $y=0\in \mathbb{S}^1$, then the support of any maximizing probability measure $\mu_\infty$ for $A(x,y) = U(x-y)$ is always contained in the set
$${\cal K}= \{ x=(...x_{-2}, x_{-1}, x_0, x_1,x_2,...)\,:\, x_i=c \in \mathbb{S}^1,  \, \, \forall i\in \mathbb{Z}\}\subset{\cal B}_i .$$

All points in ${\cal K}$ are fixed points for $\hat{\sigma}$.
The above set ${\cal K}$ can be indexed by $c\in \mathbb{S}^1$. Each fixed point $x$ in this set can be denoted by $x_c$, where $c \in \mathbb{S}^1$. The corresponding maximizing probability measure for $A$  over $(\mathbb{S}^1)^\mathbb{Z}$ is $\delta_{x_c}$.

Given any probability measure $P$ over $\mathbb{S}^1$, we can consider the probability measure $\nu$ over  $(\mathbb{S}^1)^\mathbb{Z}$ given by $\nu = \int \delta_{x_c} d P(c).$
The general maximizing probability measure for $A$ is of this form.

Suppose now that $U$ has two strict maximals at $y=0\in \mathbb{S}^1$ and at $y=\pi$. In this case, the support of any maximizing probability measure 
for $A(x,y) = U(x-y)$ is always contained in the set ${\cal K} = {\cal K}_1 \cup {\cal K}_2$, where
$${\cal K}_1= \{ x=(...x_{-2}, x_{-1}, x_0, x_1,x_2,...)\,:\, x_i=c \in \mathbb{S}^1, \, \, \forall i\in \mathbb{Z}\}\subset{\cal B}_i ,$$
and
$${\cal K}_2= \{ x=(...x_{-2}, x_{-1}, x_0, x_1,x_2,...)\,:\, x_{i+1} - x_i = \pi \in \mathbb{S}^1, \,\, \, \, \forall i\in \mathbb{Z}\}\subset{\cal B}_i .$$

The set ${\cal K}_1$ is called the set of ferromagnetic states, and,
the set ${\cal K}_2$ is called the set of anti-ferromagnetic states. The points in ${\cal K}_2$ have $\hat{\sigma}$-period equal to two. A similar result to the above is true for the general maximizing probability measure for $A$.

\bigskip

Now we will state proposition \ref{nonselec-prop} and its corollary \ref{nonselec-cor}, which, together with lemma \ref{nonselec-lemma}, will be used to prove the main result of this section, the non-selection theorem \ref{nonselec-theo}.

\noindent
\begin{proposition}\label{nonselec-prop}
 Let  $\mu_{\Lambda_n}^{\beta,\Phi}$ be the Gibbs measure in the volume $\Lambda_n$, defined by  (\ref{med-gibbs-van-enter-vol-n}).
For any fixed  $j\in\N$ and $k\in\{-n,\ldots,n-1\}$,
 if
$$
B_{k,j}=\{(\theta_{-n},\ldots,\theta_{n})\in(0,2\pi]^{2n+1}:\theta_{k}-\theta_{k+1}\in A_{j}\},
$$
then
$$
\mu_{\Lambda_n}^{\beta,\Phi}(B_{k,j})
=\frac{1}{Z(\beta)}\nu(A_{j})e^{\beta c_j},
$$
where
$$
Z(\beta)=\frac{1}{2\pi}\int_{(0,2\pi]}e^{\beta U(x)}dx,
$$
and $\nu(A_j)$ is the Lebesgue probability measure of $A_j$.
\end{proposition}

In fact, to prove theorem \ref{nonselec-theo}, we will only need to consider the Borel sets
$B_j=\{\theta_0 - \theta_1 \in A_j\} \subset{\cal B}_i$, $j\in\mathbb{N}$,
because we are interested in estimate $\mu_{\beta A}(B_j)= \int I_{B_j} \, d \mu_{\beta A}$, for each $j$,
when $\beta \to \infty$.


To state corollary \ref{nonselec-cor} we will consider the potential introduced in \cite{van}.

\begin{corollary}\label{nonselec-cor}
Let $\varepsilon >0$. Consider the special case where
$$
\tilde{U}(x)=\sum_{i=1}^{\infty}\frac{3}{2^{2i+1}}\mathds{1}_{I_{2i}}(x)
+
\sum_{i=1}^{\infty}\frac{3}{2^{2i+2}}\mathds{1}_{I_{2i+1}}(x-\pi)
+
\frac{1}{4}\mathds{1}_{I_1}\big(x-\pi\big),
$$
where $I_i=[-\frac{\varepsilon^{3^{i}}}{2},\frac{\varepsilon^{3^{i}}}{2}]$.
For each $j\in\N$ we define the ring $A_j$ as follows. If $j$ is even, then $A_{j}=A_{2i}=I_{2i}\backslash I_{2i+2}$,
and if $j$ is odd then $A_{j}=A_{2i+1}=I_{2i+1}\backslash I_{2i+3}+\pi$.
For any fixed $j\in\N$ and $k\in\{-n,\ldots,n-1\}$, we have
$$
\mu_{\beta A} (\{\theta_{k}-\theta_{k+1}\in A_{j}\})= \mu_{\Lambda_n}^{\beta,\Phi}(\{\theta_{k}-\theta_{k+1}\in A_{j}\})
=\frac{1}{Z(\beta)}\nu(A_{j})\exp\left(\frac{\beta}{2}-\frac{\beta}{2^{j+1}}\right),
$$
where
$$
Z(\beta)=\frac{1}{2\pi}\int_{(0,2\pi]}e^{\beta U(x)}dx
$$
and
$$
\nu(A_j)=\varepsilon^{3^j}-\varepsilon^{3^{j+2}}.
$$
In particular,
\begin{eqnarray}
  \nonumber \mu_{\beta A} (\{\theta_{0}-\theta_{1}\in A_{j}\}) &=& \frac{1}{Z(\beta)}\nu(A_{j})
\exp\left(\frac{\beta}{2}-\frac{\beta}{2^{j+1}}\right) \\
   &=& \frac{e^{\beta/2}}{Z(\beta)} \exp\left(-\frac{\beta}{2^{j+1}} +\log\Big(\varepsilon^{3^j}-\varepsilon^{3^{j+2}}\Big)\right)
\end{eqnarray}

\end{corollary}
\begin{figure}[h!]
 \centering
 \includegraphics[scale=0.7,keepaspectratio=true]{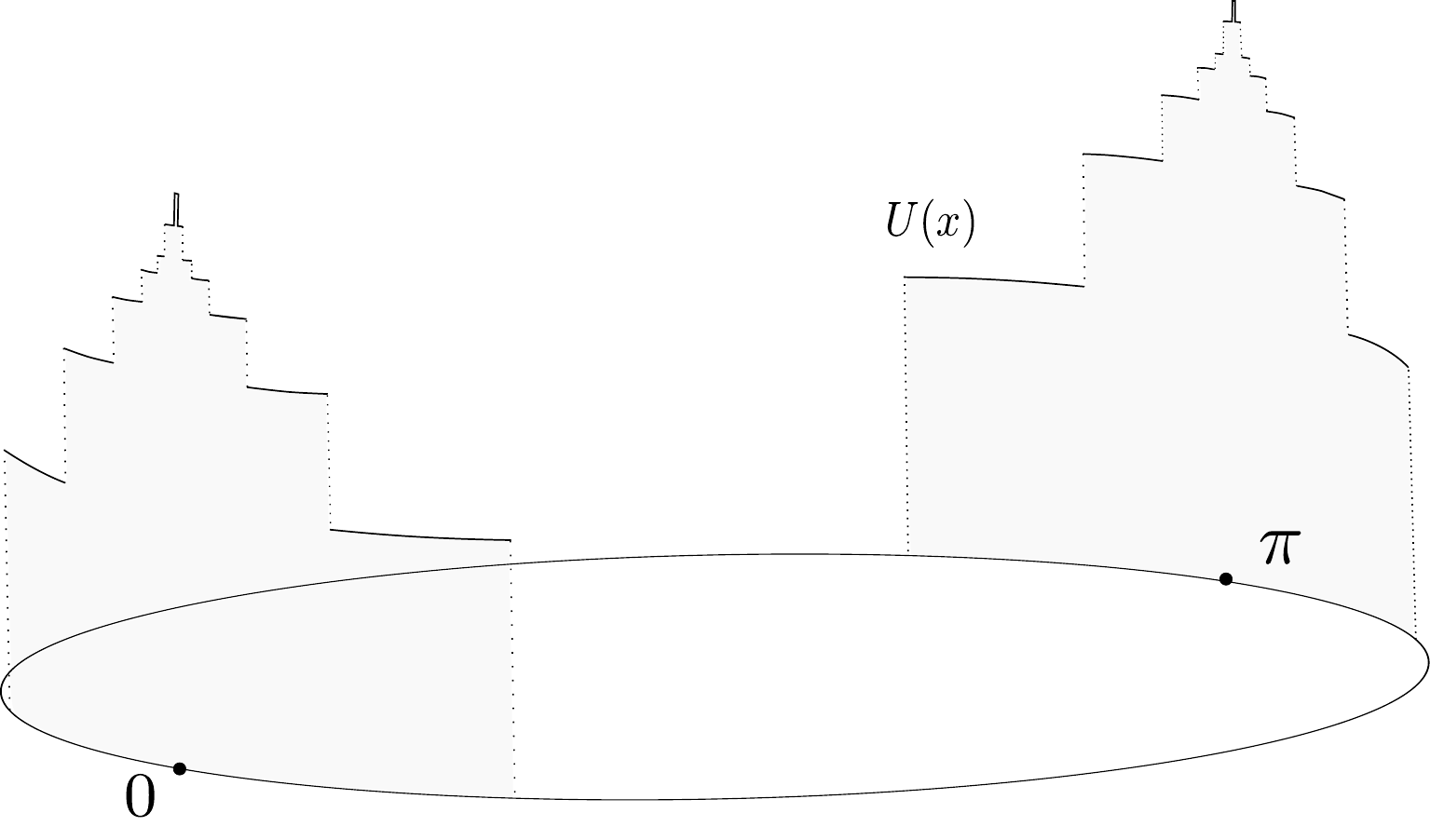}
 \caption{The graph of the potential $U$.}
 \label{fig1}
\end{figure}

\noindent{\bf Remark 6.}
Before proceeding to the proof of the proposition we remark that repeated applications of Fubini's Theorem show that
the partition function in the volume $\Lambda_n$ for the potential $A$ satisfies
$$
Z_{\Lambda_n}^{\beta,\Phi}=Z(\beta)^{2n}.
$$

\noindent {\bf Proof of Proposition.}
Let  $j\in\N$ and $k\in\{-n,\ldots,n-1\}$. 
By definition we have
$$
\mu_{\Lambda_n}^{\beta,\Phi}(\{\theta_{k}-\theta_{k+1}\in A_{j}\})=$$
$$=
\frac{1}{Z_{\Lambda_n}^{\beta,\Phi}}
\int_{(0,2\pi]^{2n+1}}
\mathds{1}_{B_{k,j}}(\theta_{-n},\ldots,\theta_{n})
\exp\left({\beta \sum_{s=-n}^{n-1}U(\theta_s-\theta_{s+1})}\right)\prod_{i=-n}^n d\nu(\theta_i).
$$

Using the properties of the exponential function, we have that the above integral is given by
\begin{equation}
\label{eq:int-inducao}
\int_{(0,2\pi]^{2n+1}}
\mathds{1}_{B_{k,j}}(\theta)
\prod_{s=-n}^{n-1}\exp\left({\beta U(\theta_s-\theta_{s+1})}\right)\prod_{i=-n}^n d\nu(\theta_i).
\end{equation}

To simplify the exposition we suppose that $k=-n$.
The following explanation can easily be modified to work in the general case just by reordering the terms, which can be done by Fubini's Theorem.
In the case $k=-n$ it follows from Fubini's Theorem that (\ref{eq:int-inducao}) is equal to
$$
\int_{(0,2\pi]^{2n}}
\mathds{1}_{B_{-n,j}}(\theta)
\prod_{s=-n}^{n-2}e^{{\beta U(\theta_{s}-\theta_{s+1})}}
\left(
\int_{(0,2\pi]}e^{{\beta U(\theta_{n-1}-\theta_{n})}}
d\nu(\theta_n)
\right)
\prod_{i=-n}^{n-1} d\nu(\theta_i).
$$
By the periodicity of $U$ it follows that the integral in parenthesis is independent of $\theta_{n-1}$ and equal to
$Z(\beta)$. Proceeding by induction, we can see that the above expression simplifies to
$$
(Z(\beta))^{2n-1}
\int_{(0,2\pi]^{2}}
\mathds{1}_{B_{-n,j}}(\theta)
e^{{\beta U(\theta_{-n}-\theta_{-n+1})}}
d\nu(\theta_{-n})d \nu(\theta_{-n+1}).
$$
To evaluate this, we consider the iterated integral where the most internal integral
is made in the variable $\theta_{-n}$, with $\theta_{-n+1}$ fixed.
For any fixed value of $\theta_{-n+1}$, whenever $\theta\in B_{-n,j}$
we have that $\theta_{-n}\in A_{j}+\theta_{-n+1}$. In this set the potential $U$ is constant, i.e.,
$$
U(\theta_{-n}-\theta_{-n+1})=c_j.
$$
From these observations the previous integral is simply
$$
(Z(\beta))^{2n-1}
\int_{(0,2\pi]}\int_{A_{j}+\theta_{-n+1}}
e^{\beta c_j}
d\nu(\theta_{-n})d \nu(\theta_{-n+1})\,,
$$
which is equal to
$$
(Z(\beta))^{2n-1}
e^{\beta c_j}
\int_{(0,2\pi]}\int_{A_{j}+\theta_{-n+1}}
d\nu(\theta_{-n})d \nu(\theta_{-n+1}).
$$
Finally by the translation invariance property of the Lebesgue measure we end up with
\begin{equation}\label{pn(beta)-simplificado}
(Z(\beta))^{2n-1}
e^{\beta c_j}
\nu(A_{j}).
\end{equation}
Dividing this value by the partition function, we get
$$
\mu_{\Lambda_n}^{\beta,\Phi}(\{\theta_{k}-\theta_{k+1}\in A_{j}\})
=\frac{\nu(A_{j})}{Z(\beta)}e^{\beta c_j}.
$$

Note that for $|k|<n$, this expression does not depend on $n$. From this follows easily that for $|k|<n$ and $j\in\N$,
$$
\mu_{\Lambda_n}^{\beta,\Phi}(\{\theta_{k}-\theta_{k+1}\in A_{j}\})=
\mu_{\beta A}(\{\theta_{k}-\theta_{k+1}\in A_{j}\}).
$$

\cqd

\vspace{0.3cm}

\noindent {\bf Proof of Corollary.}

Follows from the fact that,
if $j=2i$ and $x\in A_{2i}$ then
\begin{eqnarray*}
U(x)=\sum_{l=1}^i \frac{3}{2^{2l+1}}=\frac{3}{8}\sum_{l=0}^{i-1} \frac{1}{4^l}
&=&\frac{3}{8}\left(\frac{1-\frac{1}{4^i}}{1-\frac{1}{4}}\right)\\[0.3cm]
&=&
\frac{1}{2}\Big(1-\frac{1}{4^i}\Big)
=
\frac{1}{2}-\frac{1}{2^{2i+1}}=\frac{1}{2}-\frac{1}{2^{j+1}}.
\end{eqnarray*}
For the other hand, if $j=2i+1$ and $x\in A_{2i+1}$ we have that
\begin{eqnarray*}
U(x)=\frac{1}{4}+\sum_{l=1}^i \frac{3}{2^{2l+2}}
&=&\frac{1}{4}+\frac{3}{16}\sum_{l=0}^{i-1} \frac{1}{2^{2l}}
\\
&=&
\frac{1}{4}+\frac{3}{16}\left(\frac{1-\frac{1}{4^i}}{1-\frac{1}{4}}\right)
\\
&=&
\frac{1}{4}+\frac{1}{4}\left(1-\frac{1}{4^i}\right)
=
\frac{1}{2}-\frac{1}{2^{2i+2}}=\frac{1}{2}-\frac{1}{2^{j+1}}.
\end{eqnarray*}

\cqd

\vspace{0.3cm}

\subsection{Maximizing $\mu_{\beta A}(B_{k,j})$}\label{sec53}

Now we will present some useful calculations to compute $\mu_{\beta A}(B_{k,j})$. We point out that we will need in the future just the case $k=0.$



Let
$$
f_{\beta}(x):= -\frac{\beta}{2^{x+1}} +\log\Big(\varepsilon^{3^x}-\varepsilon^{3^{x+2}}\Big).
$$
The maximum of this function can be found by derivation with respect to $x$.
\begin{eqnarray*}
f'_{\beta}(x)&=& -\frac{d}{dx}\frac{\beta}{2^{x+1}} +\frac{d}{dx}\log\Big(\varepsilon^{3^x}-\varepsilon^{3^{x+2}}\Big)\\[0.3cm]
&=&\frac{\beta \log 2 }{2^{x+1}}+     \frac{\varepsilon^{3^x} (\log \varepsilon \log 3) 3^x-9\varepsilon^{3^{x+2}} (\log \varepsilon \log 3)3^x}{\varepsilon^{3^x}-\varepsilon^{3^{x+2}}}\\[0.3cm]
&=& \frac{\beta \log 2 }{2^{x+1}}+ (\log \varepsilon \log 3)3^x  \bigg( \frac{\varepsilon^{3^x}  -9\varepsilon^{3^{x+2}} }{\varepsilon^{3^x}-\varepsilon^{3^{x+2}}}\bigg)
\\[0.3cm]
&=& \frac{\beta \log 2 }{2^{x+1}}+ (\log \varepsilon \log 3)3^x  \bigg( \frac{\varepsilon^{3^x}  -9\varepsilon^{3^{x+2}} }{\varepsilon^{3^x}-\varepsilon^{3^{x+2}}}\bigg)\,.
\end{eqnarray*}
If $x$ is large enough the equation $f'(x)=0$ is solvable and the solution is implicitly given by
\begin{eqnarray*}
0&=&\frac{\beta \log 2 }{2^{x+1}}+ (\log \varepsilon \log 3)3^x  \bigg( \frac{\varepsilon^{3^x}  -9\varepsilon^{3^{x+2}} }{\varepsilon^{3^x}-\varepsilon^{3^{x+2}}}\bigg)\,,
\end{eqnarray*}
which is equivalent to
\begin{equation}\label{eq:betafuncaoj}
\beta=6^x\frac{-2 \log\varepsilon \log 3}{\log 2}
\bigg( \frac{\varepsilon^{3^x}-9\varepsilon^{3^{x+2}}}{\varepsilon^{3^x}  -\varepsilon^{3^{x+2}} }\bigg)\,.
\end{equation}
The fraction appearing in the above equation can be rewritten as
\begin{equation}\label{eq:theta}
\theta(\epsilon,x)\equiv\frac{\varepsilon^{3^x}}{\varepsilon^{3^x}}
\frac{1-9\varepsilon^{(3^{x+2}-3^x)}}{1-\varepsilon^{(3^{x+2}-3^x)}}
=
\frac{1-9\varepsilon^{8\cdot 3^x}}{1-\varepsilon^{8\cdot 3^x}}\,.
\end{equation}


\subsection{An important Lemma}\label{sec54}

In this subsection we present an important lemma that will help us to estimate the probability
$\mu_{\beta A}(B_{k,j})$.
\begin{lemma}\label{nonselec-lemma}
Let $(\Omega,\mathcal{B})$ be a measurable space and $(C_j)_{j\in\N}$ a measurable partition  of $\Omega$. For any positive $\beta$ let  $\mathbb{P}_{\beta}$ be a probability measure in $(\Omega,\mathcal{B})$ such that
$$\mathbb{P}_{\beta}(C_j)=\frac{1}{\bar Z(\beta)}\exp\left(-\frac{\beta}{2^{j+1}} +\log\Big(\varepsilon^{3^j}-\varepsilon^{3^{j+2}}\Big)\right)$$
where $\bar Z(\beta)$ is a normalizing constant and $\varepsilon>0$.
Given $\delta>0$, there exist an $\varepsilon_{\delta}>0$ such that, for any $0<\epsilon < \epsilon_{\delta}$,
for all $j\in\N$,  we have
$$
\mathbb{P}_{\beta_j}(C_{j})>1-\delta\,, \quad \mbox{where $\beta_j$ is given by} 
$$
$$\beta_j=6^{j}2C_{\varepsilon}\frac{\log 3}{\log 2}\theta(\varepsilon,j) \quad
with
\quad \theta(\varepsilon,j)
=
\frac{1-9\varepsilon^{8\cdot 3^j}}{1-\varepsilon^{8\cdot 3^j}} \quad
and
\quad C_{\varepsilon}=-\log\varepsilon.
$$
\end{lemma}

\noindent{\bf Remark.} A more careful analysis of the proof presented here shows that the  above lemma also works for a slightly  different potential $U$,  where we replace in the initial definition the terms $2^{j+1}$ and $3^j$,
by $(1+\delta)^{j+1}$ and $(1+\gamma)^j$, respectively, given that $0<\delta<\gamma$.


\noindent {\bf Proof.}

Note that $\theta(\varepsilon,x)$ is an increasing function of $x$, and has limit equal to $1$ when $\varepsilon \to 0$ or $x\to +\infty$.
Consider the function

 \begin{eqnarray}
  \nonumber
   f_{\beta}(x) &=& -\frac{\beta}{2^{x+1}} +\log\Big(\varepsilon^{3^x}-\varepsilon^{3^{x+2}}\Big) \\
    &=& -\frac{\beta}{2^{x+1}}-C_{\varepsilon}3^x+\log\left( 1-\varepsilon^{8.3^x}\right).\label{eq:fbeta}
 \end{eqnarray}

From \eqref{eq:betafuncaoj} and \eqref{eq:theta}
it follows that its critical point $x_0$ has to satisfy
\begin{equation}
 \label{eq:x_0-max-f}
\beta=6^{x_0}2C_{\varepsilon}\frac{\log 3}{\log 2}\theta(\varepsilon,x_0)
\end{equation}

Note that the last equation allows us to obtain the maximum  point $x_0$ of $f_{\beta}$, thus making $x_0=x_0(\beta)$ an increasing (therefore invertible) and unbounded function of $\beta$. Arguing in the inverse direction,  for each $x_0=j_0 \in \N$ we can choose $\beta=\beta(j_0)$ as the unique solution to \eqref{eq:x_0-max-f}, which means
 $j_0$ is the maximum point of $f_{\beta(j_0)}$.

Fix now $j_0\in \N$. If we set $$\kappa_{\varepsilon}^x=\log\left( 1-\varepsilon^{8.3^x}\right),$$
(note that $\kappa_{\varepsilon}^x$ is an increasing function of $x$)
it follows from \eqref{eq:fbeta}  and (\ref{eq:x_0-max-f}) that, for any $k\in\Z$
\begin{eqnarray}
f_{\beta(j_0)}(j_0+k)&=&-\frac{\beta(j_0)}{2^{j_0+k+1}}-C_{\varepsilon}3^{j_0+k}+\kappa_{\varepsilon}^{j_0+k} \\[0.4cm]
&=&-6^{j_0}2C_{\varepsilon}\frac{\log 3}{\log 2}\cdot\frac{\theta(\varepsilon,j_0)}{2^{j_0+k+1}}-C_{\varepsilon}3^{j_0+k}+\kappa_{\varepsilon}^{j_0+k}\nonumber \\[0.3cm]
\label{eq:ffacil} &=&-3^{j_0}\left[C_{\varepsilon}\frac{\log 3}{\log 2}\cdot\frac{\theta(\varepsilon,j_0)}{2^{k}}+C_{\varepsilon}3^{k}\right]+\kappa_{\varepsilon}^{j_0+k}.
\end{eqnarray}
Now we will use these identities to get an upper bound for $\frac{\mathbb{P}_{\beta(j_0)}(C_{j_0+k})}{\mathbb{P}_{\beta(j_0)}(C_{j_0})}$.
Before going to the upper bound computations we prove:

\vspace{0.1cm}
\noindent{\bf Identity 1.} For any integer $k\geq -j_0+1$, we get from \eqref{eq:ffacil} the following identity
$$
\begin{array}{lc}
\displaystyle\frac{\mathbb{P}_{\beta(j_0)}(C_{j_0+k})}{\mathbb{P}_{\beta(j_0)}(C_{j_0})}= \frac{e^{f_{\beta(j_0)}(j_0+k)}}{e^{f_{\beta(j_0)}(j_0)}}=& \\[0.6cm]
\exp\left(
-3^{j_0}\left[C_{\varepsilon}\frac{\log 3}{\log 2}\cdot\frac{\theta(\varepsilon,j_0)}{2^{k}}+C_{\varepsilon}3^{k}\right]
+3^{j_0}\left[C_{\varepsilon}\frac{\log 3}{\log2}\cdot\theta(\varepsilon,j_0)+C_{\varepsilon}\right]+\kappa_{\varepsilon}^{j_0+k}-\kappa_{\varepsilon}^{j_0}
\right)=&\\[0.3cm]
\exp\left(
-3^{j_0}C_{\varepsilon}\left[\frac{\log 3}{\log 2}\cdot\frac{\theta(\varepsilon,j_0)}{2^{k}}+3^{k}
-\frac{\log 3}{\log 2}\cdot\theta(\varepsilon,j_0)-1\right]+\kappa_{\varepsilon}^{j_0+k}-\kappa_{\varepsilon}^{j_0}
\right)=&
\\[0.3cm]
\exp\left(
-3^{j_0}C_{\varepsilon}\left[\frac{\log 3}{\log 2}\cdot\left(\frac{\theta(\varepsilon,j_0)}{2^{k}}-\theta(\varepsilon,j_0)\right)-1+3^{k}\right]+\kappa_{\varepsilon}^{j_0+k}-\kappa_{\varepsilon}^{j_0}
\right)=&
\\[0.3cm]
\exp\left(-3^{j_0}C_{\varepsilon}3^{k}\right)
\exp\left(-3^{j_0}C_{\varepsilon}\left[\frac{\log 3}{\log 2}\cdot\left(\frac{\theta(\varepsilon,j_0)}{2^{k}}-\theta(\varepsilon,j_0)\right)-1\right]+\kappa_{\varepsilon}^{j_0+k}-\kappa_{\varepsilon}^{j_0}\right)=&
\\[0.3cm]
\exp\left(-3^{j_0}C_{\varepsilon}3^{k}\right)
\exp\left(3^{j_0}C_{\varepsilon}\left[\frac{\log 3}{\log 2} \theta(\varepsilon,j_0)\cdot \left(1-\frac{1}{2^{k}}\right)+1\right]+\kappa_{\varepsilon}^{j_0+k}-\kappa_{\varepsilon}^{j_0}\right).&
\\[0.3cm]
\end{array}
$$
\\
With the above identities we are ready to show how to get the upper bounds for $\frac{\mathbb{P}_{\beta(j_0)}(C_{j_0+k})}{\mathbb{P}_{\beta(j_0)}(C_{j_0})}$. This  will be done by considering separate cases, whether $k$ is positive or negative.
\\


\noindent{\bf Case $k>0$.} In this case, using the previous identity, $\theta(\varepsilon,j_0)<1$ and $\kappa_{\varepsilon}^{j_0+k}-\kappa_{\varepsilon}^{j_0}<1$, we have
\begin{eqnarray*}
\frac{\mathbb{P}_{\beta(j_0)}(C_{j_0+k})}{\mathbb{P}_{\beta(j_0)}(C_{j_0})}
&<&
\exp\left(-3^{j_0}C_{\varepsilon}3^{k}\right)
\exp\left(3^{j_0}C_{\varepsilon}\left[\frac{\log 3}{\log 2}+1\right]+1\right)
\\[0.3cm]
&\leq&
\exp\left(-3^{j_0}C_{\varepsilon}\left[3^k-\frac{\log 3}{\log 2}-1\right]+1\right).
\end{eqnarray*}
Of course the above inequality implies, for all $k\in \N$, that
$$
\mathbb{P}_{\beta(j_0)}(C_{j_0+k})
\leq
\mathbb{P}_{\beta(j_0)}(C_{j_0})\exp\left(1-3^{j_0}C_{\varepsilon}\left[3^k-\frac{\log 3}{\log 2}-1\right]\right)
$$
and then summing over $k$ we obtain
\begin{eqnarray*}
\sum_{k=1}^{\infty}\mathbb{P}_{\beta(j_0)}(C_{j_0+k})
&\leq&
\sum_{k=1}^{\infty}
\exp\left(1-3^{j_0}C_{\varepsilon}\left[3^k-\frac{\log 3}{\log 2}-1\right]\right).
\end{eqnarray*}
In order to bound this series, 
we decompose it as follows
\begin{eqnarray*}
\exp\left(1-3^{j_0}C_{\varepsilon}
\left[2-\frac{\log 3}{\log 2}\right]\right)
+
\sum_{k=2}^{\infty}
\exp\left(1-3^{j_0}C_{\varepsilon}\left[3^k-\frac{\log 3}{\log 2}-1\right]\right).
\end{eqnarray*}
By a simple induction process one proves that $k\leq 3^k-\frac{\log 3}{\log 2}-1$, for all $k\geq 2$.
From this observation it follows the upper bound
\begin{eqnarray*}
\sum_{k=1}^{\infty}\mathbb{P}_{\beta(j_0)}(C_{j_0+k})
&\leq&
\exp\left(1-3^{j_0}C_{\varepsilon}
\left[2-\frac{\log 3}{\log 2}\right]\right)
+
\sum_{k=2}^{\infty}
\exp\left(1-3^{j_0}C_{\varepsilon}k\right)
\\[0.3cm]
&=&
\exp\left(1-3^{j_0}C_{\varepsilon}
\left[2-\frac{\log 3}{\log 2}\right]\right)
+
\frac{\exp\left(1-3^{j_0}2C_{\varepsilon}\right)}
{1-\exp\left(-3^{j_0}C_{\varepsilon}\right)}
\\[0.3cm]
&\leq&
\exp\left(1-3C_{\varepsilon}
\left[2-\frac{\log 3}{\log 2}\right]\right)
+
\frac{\exp\left(1-6C_{\varepsilon}\right)}
{1-\exp\left(-3C_{\varepsilon}\right)}.
\end{eqnarray*}
As $C_{\varepsilon}=-\log \varepsilon \to \infty$ when $\varepsilon \to 0$,
we can choose an $\varepsilon_0$ such that for all $0<\varepsilon<\varepsilon_0$, we have

\begin{equation}\label{eq:condcepsilon}
\exp\left(1-3C_{\varepsilon}
\left[2-\frac{\log 3}{\log 2}\right]\right)
+
\frac{\exp\left(1-6C_{\varepsilon}\right)}
{1-\exp\left(-3C_{\varepsilon}\right)}
<\frac{\delta}{2}.
\end{equation}
Note that $\varepsilon_0>0$ does not depend on $j_0$.

This implies that
\begin{equation}
\label{eq:est-prob-1/4-1}
\sum_{k=1}^{\infty}\mathbb{P}_{\beta(j_0)}(C_{j_0+k})<\frac{\delta}{2},
\end{equation}
for any $j_0 \in \N$, provided $0<\varepsilon<\varepsilon_0$.

\noindent
\\
\\
\noindent{\bf Case $k<0$.} From Identity 1, we have
$
\frac{\mathbb{P}_{\beta(j_0)}(C_{j_0+k})}{\mathbb{P}_{\beta(j_0)}(C_{j_0})}
$ is equal to
\begin{equation*}
\exp\left(-3^{j_0}C_{\varepsilon}3^{k}\right)
\exp\left(3^{j_0}C_{\varepsilon}\left[\frac{\log 3}{\log 2} \theta(\varepsilon,j_0)\cdot\left(1-\frac{1}{2^{k}}\right)+1\right]+\kappa_{\varepsilon}^{j_0+k}-\kappa_{\varepsilon}^{j_0}\right).
\end{equation*}
Note that we can choose
$0<\epsilon_1 \leq \epsilon_0$
such that, for all $0< \varepsilon < \epsilon_1 $ and all $j_0 \geq1$, we have
\begin{equation}\label{eq:condicepsilon-3}
\theta(\varepsilon,j_0)\frac{\log3}{\log2}-1=\frac{1-9\varepsilon^{8\cdot 3^{j_0}}}{1-\varepsilon^{8\cdot 3^{j_0}}} \frac{\log3}{\log2}-1>\frac{\frac{\log3}{\log2}-1}{2}\equiv A.
\end{equation}
As a consequence we have
\begin{equation}\label{eq:conseqcondicepsilon-3}
\theta(\varepsilon,j_0)
>\frac{\log2}{\log3}.
\end{equation} Then
$$
\frac{\log 3}{\log 2}\theta(\varepsilon,j_0)\cdot\left(1-\frac{1}{2^{k}}\right)+1<0
$$
for any $k\in\{-j_0+1,\ldots,-1\}$, and we have the following inequality, when we use  $\kappa_{\varepsilon}^{j_0+k}-\kappa_{\varepsilon}^{j_0}<0$ and $-3^{j_0}C_{\varepsilon}3^{k}<0$
\begin{equation*}\label{lim-ep-pn}
\exp\left(-3^{j_0}C_{\varepsilon}3^{k}\right)
\exp\left(3^{j_0}C_{\varepsilon}\left[\frac{\log 3}{\log 2}\theta(\varepsilon,j_0)\cdot\left(1-\frac{1}{2^{k}}\right)+1\right]
+\kappa_{\varepsilon}^{j_0+k}-\kappa_{\varepsilon}^{j_0}\right)
\leq$$
$$\leq
\exp\left(3C_{\varepsilon}\left[\frac{\log 3}{\log 2}\theta(\varepsilon,j_0)\cdot\left(1-\frac{1}{2^{k}}\right)+1\right]\right).
\end{equation*}
From this we obtain
$$
\mathbb{P}_{\beta(j_0)}(C_{j_0+k})
\leq
\exp\left(3C_{\varepsilon}\left[\frac{\log 3}{\log 2}\theta(\varepsilon,j_0)\cdot\left(1-\frac{1}{2^{k}}\right)+1\right]\right).
$$

Using this upper bound, \eqref{eq:condicepsilon-3} and \eqref{eq:conseqcondicepsilon-3} again, it follows that
\begin{eqnarray*}
\sum_{k=1}^{j_0-1}\mathbb{P}_{\beta(j_0)}(C_{j_0-k})
&\leq&
\sum_{k=1}^{\infty}
\exp\left(3C_{\varepsilon}
\left[\frac{\log 3}{\log 2}\theta(\varepsilon,j_0)\cdot\left(1-2^k\right)+1\right]\right)
\\[0.3cm]
&<&\exp\left(3C_{\varepsilon}
\left[-\frac{\log 3}{\log 2}\theta(\varepsilon,j_0)+1\right]\right)+
\sum_{k=2}^{\infty}
\exp\left(3C_{\varepsilon}
\left[\left(1-2^k\right)+1\right]\right)
\\[0.3cm]
&<&
e^{-3C_{\varepsilon}A}
+
\sum_{k=2}^{\infty}
\exp\left(3C_{\varepsilon}
\left(2-2^k\right)\right)
\\[0.3cm]
&<&
e^{-3C_{\varepsilon}A}+
\sum_{k=2}^{\infty}
\exp\left(-3C_{\varepsilon}k\right)=
e^{-3C_{\varepsilon}A}+\frac{e^{-6C_{\varepsilon}}}{1-e^{-3C_{\varepsilon}}}
\end{eqnarray*}
\\


Using again that $C_{\varepsilon}=-\log \varepsilon\to+\infty$ when $\varepsilon \to 0$, and
$A=\frac{\log3-\log2}{2\log2}>0$ we can choose $0<\varepsilon_{\delta} \leq \varepsilon_1$
such that for all $0<\varepsilon<\varepsilon_1$ we have
\begin{equation}
\label{eq:condicepsilon-2}
e^{-3\frac{\log3-\log2}{2\log2}C_{\varepsilon}}+\frac{e^{-6C_{\varepsilon}}}{1-e^{-3C_{\varepsilon}}}
<\frac{\delta}{2}\,,
\end{equation}
which implies
\begin{equation}
\label{eq:est-prob-1/4-2}
\sum_{k=1}^{j_0-1}\mathbb{P}_{\beta(j_0)}(C_{j_0-k})<\frac{\delta}{2}.
\end{equation}
Finally by (\ref{eq:est-prob-1/4-1}) and (\ref{eq:est-prob-1/4-2}) we get
$$
\sum_{k\in\N\backslash\{j_0\}}\mathbb{P}_{\beta(j_0)}(C_{k})<\delta.
$$
if $\epsilon <\epsilon_{\delta}$.
\cqd

\vspace{0.3cm}


\subsection{The non-selection theorem}

Now we are ready to state and prove the main result of this section which is due to  A. C. D. van Enter and W. M. Ruszel \cite{van}. Note that in the notation we used before the maximizing value is $m(A)= \sup U$.

\begin{theorem}\label{nonselec-theo}

For the potential $A$ described above, consider the family  of probability measures $\mu_{\beta A}$, with $\beta \in \mathbb{R}.$ Then, in the weak* topology, there is no selection of measure, that is, there is no limit for  $\mu_{\beta A}$, when $\beta \to \infty.$

\end{theorem}

\noindent{\bf Proof.}  Consider the Borel set
$$B=\{ (\theta_0- \theta_1) \in [0,\pi] \subset \mathbb{S}^1\}\subset {\cal B}_i,$$
and, the non-continuous function $I_B$. Given small $\delta$ and $\epsilon$,
we can approximate $I_B$ by a continuous function $\varphi: {\cal B}_i \to \mathbb{R}$,
where the set of points where $\varphi\neq I_B$ is contained in the small set
$${\cal D}=\{ (\theta_0- \theta_1) \in [0,\epsilon] \cup [\pi-\epsilon, \pi]\subset \mathbb{S}^1\}\subset {\cal B}_i.$$

From the above we can choose a suitable $\varphi$, and, also present two sequences $s_n$ and $t_n$, converging to infinity, such  that
$$ \int \varphi  \, d\mu_{s_n A}>1-\delta$$
and
$$ \int \varphi  \, d\mu_{t_n A}<\delta.$$

This shows that there is no limit for $\mu_{\beta A}.$

\cqd

\vspace{0.3cm}

\noindent{\bf Remark 7.} We point out that the example described above can be adapted in order to produce a continuous
potential $A$ which does not select in the limit when $\beta\to \infty$ \cite{van}.
\bigskip

\noindent


\end{document}